\documentclass[preprint,12pt,a4paper,numbers]{elsarticle}

\journal{Journal of Mathematical Analysis and Applications}
\biboptions{square,numbers,sort&compress}

\usepackage[utf8]{inputenc}
\usepackage[T1]{fontenc}
\usepackage[english]{babel}
\usepackage{amsmath,amssymb,amsthm,mathtools}
\usepackage{aliascnt}
\usepackage{enumitem}
\usepackage{booktabs,tabularx,array}
\usepackage{placeins}
\usepackage{graphicx}
\usepackage{xspace}
\usepackage{microtype}
\usepackage[hidelinks]{hyperref}
\usepackage[nameinlink,capitalise,noabbrev]{cleveref}

\hypersetup{
  breaklinks=true,
  pdftitle={Persistence and Local Extinction of Critical Branching Populations with General Lifetimes and Hunt Motions},
  pdfauthor={Ekaterina T. Kolkovska, Jose Alfredo Lopez-Mimbela, Fabio Prates Machado, Jose Hermenegildo Ramirez-Gonzalez}
}

\renewcommand{\cite}{\citet}

\newtheorem{theorem}{Theorem}[section]
\newaliascnt{proposition}{theorem}
\newtheorem{proposition}[proposition]{Proposition}
\aliascntresetthe{proposition}
\newaliascnt{lemma}{theorem}
\newtheorem{lemma}[lemma]{Lemma}
\aliascntresetthe{lemma}
\newaliascnt{corollary}{theorem}
\newtheorem{corollary}[corollary]{Corollary}
\aliascntresetthe{corollary}
\theoremstyle{remark}
\newaliascnt{remark}{theorem}
\newtheorem{remark}[remark]{Remark}
\aliascntresetthe{remark}

\crefname{theorem}{Theorem}{Theorems}
\Crefname{theorem}{Theorem}{Theorems}
\crefname{proposition}{Proposition}{Propositions}
\Crefname{proposition}{Proposition}{Propositions}
\crefname{lemma}{Lemma}{Lemmas}
\Crefname{lemma}{Lemma}{Lemmas}
\crefname{corollary}{Corollary}{Corollaries}
\Crefname{corollary}{Corollary}{Corollaries}
\crefname{remark}{Remark}{Remarks}
\Crefname{remark}{Remark}{Remarks}

\makeatletter
\@addtoreset{equation}{section}
\makeatother

\newcommand{\R}{\mathbb R}
\newcommand{\Rp}{[0,\infty)}
\newcommand{\N}{\mathbb N}
\newcommand{\E}{\mathbb E}
\newcommand{\Pp}{\mathbb P}
\newcommand{\1}{\mathbf 1}
\newcommand{\dd}{\,\mathrm d}
\newcommand{\ip}[2]{\left\langle #1,#2\right\rangle}
\newcommand{\supp}{\operatorname{supp}}
\newcommand{\weak}{\Longrightarrow}
\newcommand{\barGamma}{\overline\Gamma}
\newcommand{\cM}{\mathcal M}
\newcommand{\Hone}{\textup{(H1)}\xspace}
\newcommand{\Htwo}{\textup{(H2)}\xspace}
\newcommand{\Hthree}{\textup{(H3)}\xspace}

\begin{document}
\begin{frontmatter}

\title{Critical Branching Systems with General Lifetimes: Hunt Motions and Applications to Semilinear Equations}

\author[cimat]{Ekaterina~T. Kolkovska}
\ead{todorova@cimat.mx}
\author[cimat]{Jos\'e Alfredo L\'opez-Mimbela}
\ead{jalfredo@cimat.mx}
\author[usp]{Jos\'e Hermenegildo Ram\'irez-Gonz\'alez\corref{cor1}}
\ead{hermenegildo.ramirez@usp.br}
\cortext[cor1]{Corresponding author: Jos\'e Hermenegildo Ram\'irez-Gonz\'alez.}

\affiliation[cimat]{organization={Centro de Investigaci\'on en Matem\'aticas},
  city={Guanajuato}, country={Mexico}}
\affiliation[usp]{organization={Universidade de S\~ao Paulo}, country={Brazil}}

\begin{abstract}
We consider a critical monotype Bellman--Harris branching particle system on
$E=\mathbb{R}^d$ or $E=[0,\infty)^d$, started from a homogeneous Poisson field
with Lebesgue intensity. Particles follow a conservative Hunt process with
transition semigroup $(P_t)_{t\ge0}$, have independent lifetimes with
distribution $\Gamma$, and reproduce at death with generating function
$f(s)=s+c(1-s)^{1+\beta}$. We assume Lebesgue invariance and uniform bounds
$P_t\mathbf{1}_K=O(t^{-q})$ for compact $K$. We also impose a localization
condition: the expected number of particles in a fixed compact set at
time $t$ descended from ancestors initially outside suitable regions of
volume $O(t^q)$ can be made arbitrarily small, uniformly for $t\ge1$.
For spread-out lifetimes of finite mean, the system converges to a
nontrivial limit of Lebesgue intensity when $q\beta>1$ and becomes locally
extinct when $q\beta\le1$. For nonarithmetic lifetimes satisfying
$\overline{\Gamma}(t)\sim c_1t^{-\gamma}$, $0<\gamma\le1$, the threshold is
$q\beta=\gamma$: the limit is Poisson above it, local extinction holds
below it, and a nontrivial subsequential limit exists at equality.
We verify the hypotheses and determine recurrence--transience regimes
for symmetric and anisotropic stable processes, symmetric uniformly
elliptic diffusions, generalized squared-Bessel products on $[0,\infty)^d$,
and symmetric bilateral Bessel-type diffusions on $\mathbb{R}^d$ and
their nonsingular linear images. Finally, we derive quantitative sufficient
conditions for survival of total mass in semilinear equations with
absorption, including time-dependent diffusion rates. Here survival means
that the limiting total mass is strictly positive. For a hyperbolic
diffusion, this holds for every $\beta>0$ in the absorption term~
$-v^{1+\beta}$.
\end{abstract}

\begin{keyword}
Bellman--Harris branching \sep local extinction \sep persistence \sep
renewal equations \sep invariant measure \sep Bessel diffusion
\MSC[2020] 60J80 \sep 60K05 \sep 60J25
\end{keyword}

\end{frontmatter}

\section{Introduction}

The large-time behaviour of semilinear equations with absorption is governed
by the competition between spatial dispersal and nonlinear dissipation.
Consider
\[
 \partial_t v=\mathcal Lv-v^{1+\beta},\qquad \beta>0,
\]
where $\mathcal L$ generates a positive linear semigroup preserving the
underlying reference measure. For nonnegative solutions with finite initial
mass, absorption makes the total mass nonincreasing. The main analytic
question is whether this mass tends to zero or converges to a strictly
positive limit. We address this question through decay estimates for the
linear semigroup, obtaining quantitative lower bounds for the remaining
mass. The applications include stable, elliptic and Bessel-type generators,
as well as time-dependent diffusion rates; a hyperbolic diffusion gives a
positive limiting mass for every $\beta>0$.

The connection with branching populations can be expressed through their
Laplace functionals. For exponential lifetimes and a homogeneous Poisson
initial field, these functionals are determined by the spatial integrals
of solutions to semilinear equations with absorption. This representation
relates the long-time behaviour of the remaining mass to persistence and
local extinction of the population.

Critical spatial branching systems exhibit a competition between spatial
dispersal and genealogical clustering.  In the classical
$(d,\alpha,\beta)$ model, particles move as isotropic $\alpha$-stable
processes on $\R^d$, where $d\in\N$ and $0<\alpha\le2$, and reproduce at
death with offspring generating function $f(s)=s+c(1-s)^{1+\beta}$,
$0\le s\le1$, where $0<\beta\le1$ and $0<c<(1+\beta)^{-1}$.
Particles evolve independently, with motion, lifetime, and offspring number
independent of one another; children are born at the parent's death position.
For exponential lifetimes and an initial Poisson field of newborn ancestors
with Lebesgue intensity, the critical relation between persistence and local
extinction is $d=\alpha/\beta$; see \cite{GorostizaWakolbinger1991}.

This probabilistic question has a direct formulation in terms of a
nonlinear diffusion equation.  To describe it, let $Z_t$ be the counting
measure of descendants alive at time $t$ from one newborn ancestor at $x$,
and let $\E_x$ denote expectation for this family.  For
$\varphi\in C_c^+(\R^d)$, the space of nonnegative continuous functions
with compact support, set
\[
 a_\varphi(x):=1-e^{-\varphi(x)},\qquad
 v_\varphi(t,x):=\E_x\!\left[1-e^{-\ip{Z_t}{\varphi}}\right],
 \qquad \ip{Z_t}{\varphi}:=\int\varphi\,\dd Z_t.
\]
The function $v_\varphi(t,x)$ is denoted by $V_t^\varphi(x)$ in
\cref{sec:auxiliary}.  Write $P_t$ for the transition semigroup of the
motion, normalized so that its generator is
$\mathcal L_\alpha=-(-\Delta)^{\alpha/2}$, with Fourier symbol
$-|\xi|^\alpha$.  When lifetimes are exponential with rate
$\theta_{\mathrm e}>0$, the evolution is governed by the fractional
diffusion equation with absorption
\begin{equation}
 \partial_t v_\varphi
 =-(-\Delta)^{\alpha/2}v_\varphi
   -c\theta_{\mathrm e}v_\varphi^{1+\beta},
 \qquad v_\varphi(0,\cdot)=a_\varphi.
 \label{eq:intro-exponential-pde}
\end{equation}
For the initial data considered here, this equation is understood in its
mild form,
\begin{equation}
 v_\varphi(t,\cdot)
 =P_ta_\varphi
 -c\theta_{\mathrm e}\int_0^t
       P_{t-s}\!\left[v_\varphi(s,\cdot)^{1+\beta}\right]\dd s.
 \label{eq:intro-exponential-mild}
\end{equation}
These identities are the formulation in terms of $v_\varphi$ of the
equations in \cite[Eqs.~(2.2)--(2.4)]{GorostizaWakolbinger1991}, with the
offspring coefficient $c$ used throughout the present paper.

The relevance of the absorption term becomes explicit after spatial
integration.  With $\|u\|_1:=\int_{\R^d}|u(x)|\dd x$, define
$G_\varphi(t):=\|v_\varphi(t,\cdot)\|_1$ and
$b(\varphi):=\|a_\varphi\|_1$.  Criticality gives
$0\le v_\varphi(t,\cdot)\le P_t\varphi$, so invariance of Lebesgue measure
implies $G_\varphi(t)\le\|\varphi\|_1<\infty$.
Integrating \eqref{eq:intro-exponential-mild} then gives the mass balance
\begin{equation}
 G_\varphi(t)=b(\varphi)
 -c\theta_{\mathrm e}\int_0^t
       \|v_\varphi(s,\cdot)^{1+\beta}\|_1\dd s.
 \label{eq:intro-exponential-mass}
\end{equation}
If $X_t$ is the population obtained from independent families whose newborn
ancestors form a Poisson field of Lebesgue intensity, the Poisson formula
yields $\E e^{-\ip{X_t}{\varphi}}=e^{-G_\varphi(t)}$.  Since $G_\varphi$
is nonnegative and nonincreasing, the limiting Laplace functional is
determined, for all such test functions, by the remaining mass
\[
 \ell_\varphi:=\lim_{t\to\infty}G_\varphi(t)
 =b(\varphi)-c\theta_{\mathrm e}\int_0^\infty
       \|v_\varphi(s,\cdot)^{1+\beta}\|_1\dd s\ge0.
\]
Local extinction corresponds to $\ell_\varphi=0$ for every $\varphi$,
whereas a nontrivial limiting population leaves $\ell_\varphi>0$ for some
$\varphi$.  The analytic issue is therefore how much mass is absorbed over
time.  The mere finiteness of the absorption integral does not distinguish
these alternatives: it already follows from
\eqref{eq:intro-exponential-mass} in both regimes.

For general lifetimes, the corresponding equation retains the lifetime law
explicitly.  Let $\Gamma$ be its distribution function, with $\Gamma(0)=0$,
write $\barGamma=1-\Gamma$, and denote the associated probability measure by
$\Gamma(\dd s)$.  Conditioning on the first death gives the nonlinear
renewal equation
\begin{equation}
 \begin{aligned}
 v_\varphi(t,\cdot)
 &=\barGamma(t)P_ta_\varphi\\
 &\quad+\int_{[0,t]}P_s\!\left[
       v_\varphi(t-s,\cdot)
       -c\,v_\varphi(t-s,\cdot)^{1+\beta}\right]\Gamma(\dd s).
 \end{aligned}
 \label{eq:intro-family-renewal}
\end{equation}
This is \cite[Eq.~(10)]{VatutinWakolbinger1999} in the present notation.
For $\Gamma(t)=1-e^{-\theta_{\mathrm e}t}$ it reduces to
\eqref{eq:intro-exponential-mild}; hence the fractional absorption equation
is the exponential-lifetime case of this renewal formulation.

\cite{VatutinWakolbinger1999} considered the same branching mechanism and
stable motion with a general lifetime law.  For a nonarithmetic law of
finite mean they proved persistence for $d>\alpha/\beta$ and local
extinction for $d\le\alpha/\beta$.  For a nonarithmetic law satisfying
$\barGamma(t)\sim c_1t^{-\gamma}$, where $c_1>0$ and $0<\gamma\le1$, the
threshold becomes $d=\alpha\gamma/\beta$: above it the population converges
to the Poisson field of full Lebesgue intensity, below it the population
becomes locally extinct, and at equality a nontrivial subsequential limit
exists.  Their result at equality does not identify the limiting law.

Multitype extensions retaining symmetric stable motions include the Markovian
models with exponential lifetimes studied in
\cite{GorostizaRoellyWakolbinger1992,LopezMimbelaWakolbinger1996}, the
extinction results for multitype Bellman--Harris systems of
\cite{KeveiLopezMimbela2011}, and the finite-mean and mixed
finite--infinite mean persistence results of
\cite{KolkovskaEtAl2026,RamirezMachado2026}.  These works retain stable
spatial motions and address the additional type structure.

The present paper retains this critical reproduction mechanism and replaces
the stable motion by a conservative Hunt process on $E=\R^d$ or
$E=\Rp^d$, with transition semigroup $P_t$ and Lebesgue measure $\lambda$.
The definitions of $v_\varphi$, $G_\varphi$, and $b(\varphi)$ now refer to
families on $E$, test functions $\varphi\in C_c^+(E)$, and the norm in
$L^1(E,\lambda)$.  The same first-death decomposition gives
\eqref{eq:intro-family-renewal} with this semigroup.  If $P_t$ preserves
$\lambda$, spatial integration and resolution of the scalar renewal
equation yield
\begin{equation}
 G_\varphi(t)=b(\varphi)
 -c\int_{[0,t]}
       \|v_\varphi(t-s,\cdot)^{1+\beta}\|_1\,U(\dd s),
 \label{eq:intro-renewal-mass}
\end{equation}
where $U:=\sum_{n\ge1}\Gamma^{*n}$ is the renewal measure; $\Gamma^{*n}$ denotes
the law of the sum of $n$ independent lifetimes.  This identity is proved in
\cref{lem:Laplace}.  For exponential lifetimes,
$U(\dd s)=\theta_{\mathrm e}\dd s$, recovering
\eqref{eq:intro-exponential-mass}.  In general, the Laplace functionals are
therefore determined by the renewal convolution of the $L^1$ norm of
$v_\varphi^{1+\beta}$.

The assumptions on $P_t$ are stated in \Hone--\Hthree.  They require
invariance of Lebesgue measure, a uniform estimate
$P_t\1_K\le C_K(1+t)^{-q}$ for each compact $K\subset E$, with an exponent
$q>0$ and a constant $C_K<\infty$, and a localization property on regions
of volume $O(t^q)$.  These are the properties of the stable
semigroup that enter the renewal argument of
\cite{VatutinWakolbinger1999}.  They control the nonlinear term in
\eqref{eq:intro-renewal-mass} through dispersal and localization estimates;
their distinct roles and logical independence are established in
\cref{sec:model}.

Under these assumptions the stable exponent $d/\alpha$ is replaced by $q$.
For spread-out lifetimes of finite mean, the population persists when
$q\beta>1$ and becomes locally extinct when $q\beta\le1$.  In this regime,
the spread-out assumption is stronger than the nonarithmetic hypothesis in
\cite{VatutinWakolbinger1999}; it is used below through Stone's decomposition
of the renewal measure.  For a nonarithmetic lifetime law with
$\barGamma(t)\sim c_1t^{-\gamma}$, the threshold is $q\beta=\gamma$.
The renewal estimates needed at the endpoint $\gamma=1$ are included.

The semigroup hypotheses hold for several standard and nonhomogeneous Markov
motions.  Besides full-dimensional symmetric stable processes, they include
anisotropic stable products with different indices and all their nonsingular
linear images; in particular, the Brownian case allows every positive
definite covariance matrix.  They also hold for symmetric uniformly elliptic
divergence-form diffusions.  A further class is obtained by a power
transformation of squared-Bessel processes, parametrized by $\rho>0$ as in
\cref{sec:Bessel}.  The resulting one-dimensional
kernel is symmetric with respect to Lebesgue measure and decays locally as
$t^{-1/\rho}$, so the independent $d$-fold product has exponent $q=d/\rho$.
When $\rho>1$, the symmetric extension through zero gives a conservative
Bessel-type diffusion on the whole line, and nonsingular linear images yield
examples on $\R^d$ with dependent coordinates.

For each listed example we also determine recurrence and transience of the
individual motion.  The estimate \Htwo gives transience whenever $q>1$.
For the symmetric classes in \cref{tab:examples}, recurrence holds when
$q\le1$, while the Brownian motion with nonzero drift considered after the
table shows that this conclusion is not implied by \Hone--\Hthree.  This
comparison is separate from the branching threshold: finite-mean persistence
implies transience of the individual motion, while for heavy-tailed lifetimes
a Poisson limit may occur in a recurrent regime.

The paper is organized as follows.  \Cref{sec:model} introduces the model,
states the hypotheses on the motion, and gives the extinction--survival
results.  \Cref{sec:examples} first verifies the motion hypotheses and then
treats recurrence and transience.  The estimates for one family and the
renewal equations are proved in \cref{sec:auxiliary}; the proofs of the
extinction--survival results are given in \cref{sec:proofs}.
Finally, \cref{sec:semilinear-survival} applies the semigroup decay estimates
to persistence of the total mass of semilinear equations with absorption,
including equations with time-dependent diffusion rates.

\section{Model and extinction--survival results}
\label{sec:model}

Fix $d\in\N$.  Let $E=\R^d$ or $E=\Rp^d$, with Borel
$\sigma$-field and Lebesgue measure $\lambda$.  Let $W=(W_t)_{t\ge0}$ be
a conservative Hunt process with Markov family $(\Pp_x)_{x\in E}$ and
semigroup $P_th(x)=\E_x[h(W_t)]$, $h\ge0$.  Each particle moves as a copy
of $W$ from its birth position, has lifetime distribution $\Gamma$ with
$\Gamma(0)=0$, and at death is replaced at its death position by
$\zeta\in\N_0$ children.  For each particle, its motion, lifetime, and
offspring number are mutually independent; the corresponding marks of
distinct particles are independent.  Write
$\barGamma(t)=1-\Gamma(t)$ and $\Gamma(\dd t)$ for the corresponding
probability measure.  The offspring generating function is
$f(s)=\E[s^\zeta]=s+c(1-s)^{1+\beta}$, where $0<\beta\le1$ and
$0<c<(1+\beta)^{-1}$.  This is exactly the nondegenerate critical mechanism
of \cite[Eq.~(1)]{VatutinWakolbinger1999}.  Since $f(s)>s$ for $0\le s<1$,
the associated Galton--Watson process becomes extinct after a finite
generation almost surely.  We use the standard spatial Bellman--Harris model
of \cite[Section~2]{VatutinWakolbinger1999}: particles reproduce only at
death, all children are born at the parent's death position, and distinct
families evolve independently.  

Let $Z_t^x$ be the counting measure of descendants alive at time $t$ from one
newborn ancestor at $x$.  If $\Pi$ is a Poisson random measure of intensity
$\lambda$, attach independent one-family marks and set
$X_t=\sum_{x\in\Pi}Z_t^x$.
Let $\cM(E)$ be the space of nonnegative Radon measures on $E$, equipped
with the vague topology.  Every convergence statement for $X_t$ below means
convergence in distribution on this space.  Throughout, local extinction
means $X_t\weak0$.  In particular,
$\Pp\{X_t(K)>0\}\to0$ for every compact $K$: choose
$\psi\in C_c^+(E)$ with $\psi\ge\1_K$ and use
$\{X_t(K)>0\}\subset\{\ip{X_t}{\psi}\ge1\}$.
Recall that $\Gamma$ is spread-out if some $\Gamma^{*n}$ has a nonzero
absolutely continuous component, and nonarithmetic if it is not supported by
$h\N_0$ for any $h>0$.

\paragraph{Hypotheses}
For some $q>0$, the semigroup satisfies
\begin{enumerate}[label=\textup{(H\arabic*)},leftmargin=3.3em,itemsep=3pt,topsep=3pt]
\item \emph{Invariant reference measure:}
$\int_E P_th\,\dd\lambda=\int_Eh\,\dd\lambda$ for every $t\ge0$ and every
nonnegative Borel $h$, with either side allowed to be infinite.
\item \emph{Uniform local dispersal:} for each compact $K\subset E$, there is
$0<C_K<\infty$ such that $\sup_xP_t\1_K(x)\le C_K(1+t)^{-q}$, $t\ge0$.
\item \emph{Dispersal-scale localization:} there are Borel sets
$D_t(N)\subset E$, $t,N\ge1$, and constants $C_N<\infty$ such that
$\lambda(D_t(N))\le C_Nt^q$ and, for every
compact $K\subset E$,
$\lim_{N\to\infty}\sup_{t\ge1}\int_{E\setminus D_t(N)}P_t\1_K\,\dd\lambda=0$.
\end{enumerate}

\paragraph{Distinct roles and logical independence of
\textup{(H1)}--\textup{(H3)}}
The three semigroup hypotheses have distinct roles in the argument.  Condition
\Hone preserves Lebesgue intensity under the motion and permits the spatially
integrated equations to be reduced to scalar renewal equations.  Condition
\Htwo gives the uniform local decay that determines the dispersal exponent
$q$, whereas \Hthree ensures that, uniformly in time, the initial positions
that may contribute to a fixed compact set can be confined to sets of volume
of order $t^q$.

We next record that none of these assumptions is a consequence of the other
two.  Since the statement of the hypotheses requires the existence of a
single exponent $q>0$ for which both \Htwo and \Hthree hold, the exponent must
be quantified carefully.  For a conservative Hunt semigroup $P$, define
\[
\begin{aligned}
 \mathcal Q_2(P)
 &:=\{q>0:P\text{ satisfies \Htwo with exponent }q\},\\
 \mathcal Q_3(P)
 &:=\{q>0:P\text{ satisfies \Hthree with exponent }q\}.
\end{aligned}
\]
These sets have opposite monotonicity properties:
\[
 q_0\in\mathcal Q_2(P)\ \Longrightarrow\
 (0,q_0]\subseteq\mathcal Q_2(P),
 \qquad
 q_0\in\mathcal Q_3(P)\ \Longrightarrow\
 [q_0,\infty)\subseteq\mathcal Q_3(P).
\]
Indeed, the first implication follows from
$(1+t)^{-q_0}\le(1+t)^{-q}$ for $0<q\le q_0$, and the second follows by
using the same localization sets and observing that $t^{q_0}\le t^q$ for
$t\ge1$ and $q\ge q_0$.  Thus it would not be sufficient merely to exhibit
one exponent for which one condition holds and the other fails.  The
examples below prove the stronger assertion that changing the exponent
cannot restore the missing hypothesis.

\begin{enumerate}[
 label=\textup{(\roman*)},
 leftmargin=2.2em,
 itemsep=5pt,
 topsep=4pt
]

\item \emph{Conditions \Hone and \Hthree do not imply \Htwo.}
On $E=\R^d$, consider the constant process $W_t^x=x$.  It is deterministic,
has continuous paths and infinite lifetime, and satisfies the strong Markov
property.  Hence it is a conservative Hunt process.  Its semigroup is
$P_th=h$, and therefore \Hone holds.  For an arbitrary $q>0$, let
\[
 D_t(N)=B\bigl(0,Nt^{q/d}\bigr),\qquad t,N\ge1.
\]
Then
\[
 \lambda(D_t(N))=\lambda(B(0,1))N^dt^q.
\]
If $K\subset B(0,R)$ is compact, then, for every $N\ge R$ and $t\ge1$,
\[
 \int_{\R^d\setminus D_t(N)}P_t\1_K\,\dd\lambda
 =\int_{\R^d\setminus D_t(N)}\1_K\,\dd\lambda=0.
\]
Consequently, \Hthree holds for every $q>0$.  On the other hand, for every
nonempty compact set $K$,
\[
 \sup_{x\in\R^d}P_t\1_K(x)=1,\qquad t\ge0,
\]
so \Htwo fails for every $q>0$.  Therefore
\[
 \mathcal Q_3(P)=(0,\infty),
 \qquad
 \mathcal Q_2(P)=\varnothing.
\]

\item \emph{Conditions \Hone and \Htwo do not imply \Hthree.}
On $E=\R^2$, let
\[
 A=\begin{pmatrix}1&0\\0&-1\end{pmatrix}
\]
and consider the linear diffusion
\[
 \dd W_t=AW_t\,\dd t+\dd B_t,
\]
where $B$ is a standard Brownian motion in $\R^2$.  The coefficients are
globally Lipschitz and have at most linear growth, so the equation has a
unique nonexplosive strong solution with continuous paths and a Feller
transition semigroup.  Hence $W$ is a conservative continuous Feller
process, and therefore a Hunt diffusion.  Starting from $x$, it has the
representation
\[
 W_t^x=e^{tA}x+Z_t,
 \qquad
 Z_t:=\int_0^t e^{(t-s)A}\,\dd B_s,
\]
where $Z_t$ is centered Gaussian with covariance matrix
\[
 Q_t=
 \begin{pmatrix}
  (e^{2t}-1)/2&0\\
  0&(1-e^{-2t})/2
 \end{pmatrix},
 \qquad
 \det Q_t=\sinh^2 t.
\]
Since $\det(e^{tA})=e^{t\operatorname{tr}A}=1$, for every nonnegative Borel
function $h$, Tonelli's theorem and the change of variables
$y=e^{tA}x+Z_t$ give
\[
 \int_{\R^2}P_th(x)\,\dd x
 =
 \E\!\left[
  \int_{\R^2}h(e^{tA}x+Z_t)\,\dd x
 \right]
 =
 \int_{\R^2}h(y)\,\dd y.
\]
Thus \Hone holds.  The transition density is
\[
 p_t(x,y)
 =
 \frac{1}{2\pi(\det Q_t)^{1/2}}
 \exp\!\left[
  -\frac12
  \bigl(y-e^{tA}x\bigr)^\top Q_t^{-1}
  \bigl(y-e^{tA}x\bigr)
 \right],
\]
and hence
\[
 \sup_{x,y\in\R^2}p_t(x,y)
 =
 \frac{1}{2\pi\sinh t},
 \qquad t>0.
\]
It follows that, for every compact $K\subset\R^2$ and $t\ge1$,
\[
 \sup_{x\in\R^2}P_t\1_K(x)
 \le
 \frac{\lambda(K)}{2\pi\sinh t}.
\]
Since $(\sinh t)^{-1}$ decreases faster than every negative power of $t$,
and since $P_t\1_K\le1$ for $0\le t\le1$, condition \Htwo holds for every
$q>0$.

We claim that \Hthree fails for every $q>0$.  Fix $q>0$ and suppose, to the
contrary, that there are sets $D_t(N)$ satisfying \Hthree with this exponent.
Choose a compact set $K$ with $0<\lambda(K)<\infty$.  For each fixed $N$,
\[
 \begin{aligned}
 \int_{D_t(N)}P_t\1_K(x)\,\dd x
 &\le
 \lambda(D_t(N))
 \sup_{x\in\R^2}P_t\1_K(x)\\
 &\le
 \frac{C_N\lambda(K)t^q}{2\pi\sinh t}
 \longrightarrow0,
 \qquad t\to\infty.
 \end{aligned}
\]
By \Hone,
\[
 \int_{\R^2}P_t\1_K\,\dd\lambda=\lambda(K),
\]
and therefore, for every fixed $N$,
\[
 \int_{\R^2\setminus D_t(N)}P_t\1_K\,\dd\lambda
 \longrightarrow\lambda(K),
 \qquad t\to\infty.
\]

Consequently, for every fixed $N$,
\[
 \lim_{t\to\infty}
 \int_{\R^2\setminus D_t(N)}
 P_t\1_K\,\dd\lambda
 =
 \lambda(K).
\]
Moreover, since
\[
 0\le
 \int_{\R^2\setminus D_t(N)}
 P_t\1_K\,\dd\lambda
 \le
 \int_{\R^2}P_t\1_K\,\dd\lambda
 =
 \lambda(K),
\]
we obtain
\[
 \sup_{t\ge1}
 \int_{\R^2\setminus D_t(N)}
 P_t\1_K\,\dd\lambda
 =
 \lambda(K)
\]
for every $N$.  Therefore,
\[
 \lim_{N\to\infty}
 \sup_{t\ge1}
 \int_{\R^2\setminus D_t(N)}
 P_t\1_K\,\dd\lambda
 =
 \lambda(K)>0,
\]
which contradicts \Hthree.  Hence
\[
 \mathcal Q_2(P)=(0,\infty),
 \qquad
 \mathcal Q_3(P)=\varnothing.
\]

\item \emph{Conditions \Htwo and \Hthree do not imply \Hone.}
Let $B$ be a standard Brownian motion on $\R$, fix $0<\varepsilon<1$, and set
\[
 T(x)=x+\varepsilon\sin x,
 \qquad
 m:=1-\varepsilon,
 \qquad
 M:=1+\varepsilon.
\]
Then $T$ is an increasing $C^1$-diffeomorphism of $\R$ and
\[
 m\le T'(x)\le M,
 \qquad x\in\R.
\]
Define the homeomorphic image of Brownian motion by
\[
 \widetilde W_t^x
 =
 T\bigl(T^{-1}(x)+B_t\bigr),
\]
whose semigroup is
\[
 \widetilde P_t f(x)
 =
 P_t^B(f\circ T)(T^{-1}x),
\]
where $(P_t^B)_{t\ge0}$ denotes the standard Brownian semigroup.  Since $T$
is a homeomorphism, the map $Uf=f\circ T$ is an isometric isomorphism of
$C_0(\R)$, and $\widetilde P_t=U^{-1}P_t^BU$ is a conservative Feller
semigroup.  Moreover, $\widetilde W$ has continuous paths.  Hence it is a
conservative continuous Feller process, and therefore a Hunt diffusion.  If
$K\subset\R$ is compact and $K^*:=T^{-1}(K)$, then, for $t\ge1$,
\[
 \sup_{x\in\R}\widetilde P_t\1_K(x)
 =
 \sup_{z\in\R}P_t^B\1_{K^*}(z)
 \le
 \frac{\lambda(K^*)}{\sqrt{2\pi t}}
 \le
 \frac{\lambda(K)}{m\sqrt{2\pi t}}.
\]
Together with the bound $\widetilde P_t\1_K\le1$ on $0\le t\le1$, this
proves \Htwo with $q=1/2$.

For $t,N\ge1$, set
\[
 D_t(N)=T([-N\sqrt t,N\sqrt t]).
\]
Then
\[
 \lambda(D_t(N))\le2MN\sqrt t.
\]
Let $K^*\subset[-R,R]$.  By the change of variables $x=T(z)$, the bound
$T'\le M$, Tonelli's theorem, and symmetry of the Brownian transition
density, for $N>R$ and $t\ge1$ we obtain
\[
 \begin{aligned}
 \int_{\R\setminus D_t(N)}
 \widetilde P_t\1_K(x)\,\dd x
 &\le
 M\int_{\R\setminus[-N\sqrt t,N\sqrt t]}
 P_t^B\1_{K^*}(z)\,\dd z\\
 &=
 M\int_{K^*}
 \Pp\{y+B_t\notin[-N\sqrt t,N\sqrt t]\}\,\dd y\\
 &\le
 M\lambda(K^*)\Pp\{|B_1|>N-R\}.
 \end{aligned}
\]
The last expression tends to zero as $N\to\infty$, uniformly in $t\ge1$.
Thus \Hthree also holds with the same exponent $q=1/2$.

It remains to show that Lebesgue measure is not invariant.  For every
nonnegative compactly supported Borel function $f$, symmetry of the Brownian
semigroup and the change of variables $x=T(z)$ give
\[
 \begin{aligned}
 \int_{\R}\widetilde P_t f(x)\,\dd x
 &=
 \int_{\R}P_t^B(f\circ T)(z)T'(z)\,\dd z\\
 &=
 \int_{\R}(f\circ T)(z)P_t^BT'(z)\,\dd z\\
 &=
 \int_{\R}f(y)
 \frac{P_t^BT'(T^{-1}y)}
      {T'(T^{-1}y)}
 \,\dd y.
 \end{aligned}
\]
Therefore, invariance of Lebesgue measure would imply
$P_t^BT'=T'$ almost everywhere.  However,
\[
 T'(z)=1+\varepsilon\cos z,
 \qquad
 P_t^BT'(z)
 =
 1+\varepsilon e^{-t/2}\cos z,
\]
which are not equal for any $t>0$.  Hence \Hone fails.

\end{enumerate}

The three examples prove, within the class of conservative Hunt processes,
that none of \Hone--\Hthree follows from the other two, even when the
exponent $q$ is allowed to vary.

\begin{theorem}\label{thm:finite}
Assume \Hone--\Hthree.  Suppose that $\Gamma$ is spread-out and has
finite mean $\mu=\E\sigma=\int_0^\infty\barGamma(t)\,\dd t$.  Then:
\begin{enumerate}[label=\textup{(\alph*)},leftmargin=2em,itemsep=2pt,topsep=2pt]
\item if $q\beta>1$, then $X_t\weak X_\infty$ for a nontrivial random Radon
measure satisfying $\E\ip{X_\infty}{\varphi}=\int_E\varphi\,\dd\lambda$ for
all $\varphi\in C_c^+(E)$;
\item if $q\beta\le1$, then $X_t\weak0$.
\end{enumerate}
\end{theorem}

\begin{theorem}\label{thm:infinite}
Assume \Hone--\Hthree.  Suppose that $\Gamma$ is nonarithmetic and
\begin{equation}
 \barGamma(t)\sim c_1t^{-\gamma},
 \qquad c_1>0,\quad 0<\gamma\le1.
 \label{eq:tail}
\end{equation}
Then:
\begin{enumerate}[label=\textup{(\alph*)},leftmargin=2em,itemsep=2pt,topsep=2pt]
\item if $q\beta>\gamma$, then $X_t\weak\Pi_\lambda$, a Poisson random
measure with intensity $\lambda$;
\item if $q\beta<\gamma$, then $X_t\weak0$;
\item if $q\beta=\gamma$, there are $t_n\to\infty$ and a random Radon
measure $X_\infty$ such that $X_{t_n}\weak X_\infty$ and
$\Pp\{X_\infty\ne0\}>0$.
\end{enumerate}
\end{theorem}

\section{Markov process examples}
\label{sec:examples}

We first verify \Hone--\Hthree for four classes of conservative Hunt
processes and identify the corresponding exponent \(q\).  Recurrence and
transience are treated separately in \cref{subsec:recurrence-transience}.
The estimate \Htwo gives a general sufficient condition for transience when
\(q>1\); recurrence is then established separately for each symmetric
class.  Throughout this section, \(GL_d(\R)\) denotes the group of
invertible \(d\times d\) real matrices.

For $F=\R^d$ or $\Rp^d$, $C_0(F)$ is the space of continuous functions
vanishing at infinity, equipped with the supremum norm
$\|h\|_\infty=\sup_{x\in F}|h(x)|$.  We write $C_c(F)$ for the continuous
functions with compact support and $C_c^+(F)$ for its nonnegative cone.
For an open set $O$, $C_c^k(O)$ consists of $k$ times continuously
differentiable functions with compact support contained in $O$.
In $C_c^2([0,\infty))$, the derivatives extend continuously to zero from
the right.

\paragraph{Two recurring facts.}
If a conservative semigroup admits a transition density satisfying
$p_t(x,y)=p_t(y,x)$, then Lebesgue measure is invariant, the semigroup is
reversible, and its $L^2$-extension is self-adjoint.  We shall also use the
localization identity
\[
 \int_{E\setminus D}P_t\1_K(x)\lambda(\dd x)
 =\int_K\Pp_y\{W_t\notin D\}\lambda(\dd y),
\]
valid for Borel $D\subset E$ and compact $K\subset E$.  Both statements
follow directly from Tonelli's theorem, symmetry, and conservativity.
Moreover, once an estimate of the form
$\sup_xP_t\1_K(x)\le C_Kt^{-q}$ has been proved for $t\ge1$, the bound
$P_t\1_K\le1$ extends it to $C_K'(1+t)^{-q}$ for all $t\ge0$.  We use this
extension without further comment.

\begin{lemma}\label{lem:linear-images}
Let $W$ be a conservative Hunt process on $\R^d$ satisfying
\Hone--\Hthree with exponent $q$ and localization sets $D_t(N)$.  For
$A\in GL_d(\R)$, the process $W_t^A=AW_t$ also satisfies
\Hone--\Hthree with exponent $q$ and localization sets $AD_t(N)$.  If the
semigroup of $W$ is symmetric with respect to Lebesgue measure, then so is
the semigroup of $W^A$.  If $U$ and $U^A$ denote the corresponding
potential operators, then
\begin{equation}
 U^A\1_K(Ax)=U\1_{A^{-1}K}(x)
 \label{eq:linear-image-potential}
\end{equation}
for every Borel set $K\subset\R^d$.
\end{lemma}

\begin{proof}
Since $A$ is a homeomorphism, $W^A$ is again a conservative Hunt process,
with semigroup
\[
 P_t^Ah(x)=P_t(h\circ A)(A^{-1}x).
\]
A change of variables shows that \Hone is preserved.  For every compact
$K\subset\R^d$,
\[
 \begin{aligned}
 \sup_xP_t^A\1_K(x)&=\sup_zP_t\1_{A^{-1}K}(z),\\
 \int_{(AD_t(N))^c}P_t^A\1_K\,\dd\lambda
 &=|\det A|\int_{D_t(N)^c}P_t\1_{A^{-1}K}\,\dd\lambda,
 \end{aligned}
\]
and $\lambda(AD_t(N))=|\det A|\lambda(D_t(N))$.  Hence \Htwo and
\Hthree hold with the same exponent.  The same change of variables in the
bilinear symmetry identity proves symmetry of the transformed semigroup,
and integration of the semigroup identity over time gives
\eqref{eq:linear-image-potential}.
\end{proof}

\subsection{Symmetric stable motions and anisotropic stable products}

We first record two classes of stable motions satisfying \Hone--\Hthree.
Although \cite{VatutinWakolbinger1999} refers simply to a symmetric
\(\alpha\)-stable motion, the estimates used there are those for the
spherically symmetric stable process with generator
\(-(-\Delta)^{\alpha/2}\), corresponding to the isotropic case. The first
class below removes this isotropy restriction and consists of arbitrary
full-dimensional symmetric stable L\'evy processes, while the second is
obtained from products of one-dimensional stable motions by an invertible
linear transformation.

\begin{proposition}\label{prop:stable}
\begin{enumerate}[label=\textup{(\alph*)},leftmargin=2em,itemsep=2pt,topsep=2pt]

\item Let \(0<\alpha\le2\). A full-dimensional symmetric
\(\alpha\)-stable L\'evy process on \(\R^d\), that is, one whose law at
time one is not supported on a proper linear subspace of \(\R^d\),
satisfies \Hone--\Hthree with \(q=d/\alpha\) and
\(D_t(N)=B(0,Nt^{1/\alpha})\).

\item Let \(S^{(i)}\), \(i=1,\ldots,d\), be independent nondegenerate
symmetric \(\alpha_i\)-stable L\'evy processes on \(\R\), with
\(0<\alpha_i\le2\), and let \(A\in GL_d(\R)\). Then
\(S_t^A=A(S_t^{(1)},\ldots,S_t^{(d)})\) satisfies \Hone--\Hthree with
\(q=\sum_{i=1}^d\alpha_i^{-1}\) and
\[
D_t(N)=
A\left(\prod_{i=1}^d[-Nt^{1/\alpha_i},Nt^{1/\alpha_i}]\right).
\]
In particular, this class contains every centered Brownian motion on
\(\R^d\) with positive definite covariance matrix.

\end{enumerate}
\end{proposition}

\begin{proof}

Let \(W\) be the process in part~\textup{(a)} and write
\[
\E e^{i\xi\cdot W_t}=e^{-t\Psi(\xi)},\qquad \xi\in\R^d.
\]
Suppose first that \(0<\alpha<2\). By
\cite[Theorem~14.13]{Sato1999}, there exists a unique finite symmetric
measure \(\sigma\) on \(\mathbb S^{d-1}\) such that
\[
\Psi(\xi)=\int_{\mathbb S^{d-1}}
|\langle\xi,\theta\rangle|^\alpha\,\sigma(d\theta).
\]
The finiteness of \(\sigma\) implies, by dominated convergence, that
\(\Psi\) is continuous, and its representation gives
\(\Psi(r\xi)=|r|^\alpha\Psi(\xi)\). If \(\Psi(\xi_0)=0\) for some
\(\xi_0\ne0\), this homogeneity would give
\(\E e^{ir\langle\xi_0,W_1\rangle}=1\) for every \(r\in\R\).
Then \(\langle\xi_0,W_1\rangle=0\) almost surely, contradicting full
dimensionality. Hence \(\Psi\) is strictly positive on
\(\mathbb S^{d-1}\), and compactness and homogeneity yield
\[
a|\xi|^\alpha\le\Psi(\xi)\le b|\xi|^\alpha,
\qquad \xi\in\R^d,
\]
for some \(0<a\le b<\infty\).
For \(\alpha=2\), \(W_1\) is centered Gaussian
\cite[Theorem~14.1]{Sato1999}, with a positive definite covariance
matrix \(\Sigma\) by full dimensionality. The same homogeneity and
bounds follow from
\[
\Psi(\xi)=\tfrac12\langle\Sigma\xi,\xi\rangle,
\qquad
a=\tfrac12\lambda_{\min}(\Sigma),\quad
b=\tfrac12\lambda_{\max}(\Sigma).
\]

Since \(\Psi(\xi)\ge a|\xi|^\alpha\), for every \(t>0\),
\[
\int_{\R^d}e^{-t\Psi(\xi)}\,d\xi
\le
\int_{\R^d}e^{-at|\xi|^\alpha}\,d\xi<\infty.
\]
Thus the characteristic function of \(W_t\) is integrable. By
\cite[Proposition~2.5(xii)]{Sato1999}, the law of \(W_t\) has a bounded
continuous density
\[
p_t(z)=(2\pi)^{-d}\int_{\R^d}
e^{-i\xi\cdot z}e^{-t\Psi(\xi)}\,d\xi.
\]
Changing variables \(\eta=t^{1/\alpha}\xi\) in this absolutely
convergent integral and using homogeneity gives
\[
 p_t(z)=t^{-d/\alpha}p_1(t^{-1/\alpha}z).
\]
The symmetry of \(\Psi\) also gives \(p_t(z)=p_t(-z)\).
Every L\'evy process is a Hunt process
\cite[Theorem~3.6.11]{Applebaum2009}, and its semigroup is conservative.
The convolution form of the semigroup and translation invariance of
Lebesgue measure give \Hone.  Moreover, for every compact
\(K\subset\R^d\) and \(t\ge1\),
\[
\sup_{x\in\R^d}P_t\mathbf1_K(x)
\le
\lambda(K)\|p_t\|_\infty
=
\lambda(K)\|p_1\|_\infty t^{-d/\alpha}.
\]

Thus \Htwo holds with \(q=d/\alpha\).

For \Hthree, take \(K\subset B(0,R)\) and
\(D_t(N)=B(0,Nt^{1/\alpha})\). By the localization identity at the
beginning of the section and \(W_t\overset{d}=t^{1/\alpha}W_1\),
for \(t\ge1\) and \(N>R\),
\[
\begin{aligned}
\int_{D_t(N)^c}P_t\mathbf1_K(x)\,dx
&=\int_K\Pp\{y+t^{1/\alpha}W_1\notin D_t(N)\}\,dy\\
&\le\lambda(K)\Pp\{|W_1|>N-R\}.
\end{aligned}
\]
Indeed, the event in the first line implies
\(|W_1|>N-|y|t^{-1/\alpha}\ge N-R\).
The bound tends to zero as \(N\to\infty\), uniformly in \(t\ge1\),
and \(\lambda(D_t(N))=\lambda(B(0,N))t^{d/\alpha}\). This proves
\Hthree.

For part~\textup{(b)}, set
\(S_t=(S_t^{(1)},\ldots,S_t^{(d)})\) and
\(q=\sum_{i=1}^d\alpha_i^{-1}\). The product \(S\) is again a
conservative L\'evy process, hence Hunt by the result used in
part~\textup{(a)}. Applying that part in dimension one and using
independence gives the symmetric density
\[
p_t(z)=t^{-q}\prod_{i=1}^d
p_1^{(i)}\!\left(t^{-1/\alpha_i}z_i\right),
\qquad \|p_1^{(i)}\|_\infty<\infty.
\]
Its convolution semigroup satisfies \Hone, and
\[
\sup_xP_t\mathbf1_K(x)
\le\lambda(K)t^{-q}\prod_{i=1}^d\|p_1^{(i)}\|_\infty,
\qquad t\ge1,
\]
proves \Htwo.

Set
\[
\begin{aligned}
Q_t(N)&=\prod_{i=1}^d[-Nt^{1/\alpha_i},Nt^{1/\alpha_i}]
      =\Delta_t[-N,N]^d,
\\
\Delta_t&=\operatorname{diag}
(t^{1/\alpha_1},\ldots,t^{1/\alpha_d}).
\end{aligned}
\]
Then \(\lambda(Q_t(N))=(2N)^dt^q\). For
\(K\subset[-R,R]^d\), the localization identity and
\(S_t\overset{d}=\Delta_tS_1\) give
\[
\begin{aligned}
\int_{Q_t(N)^c}P_t\mathbf1_K(x)\,dx
&=\int_K\Pp\{\Delta_t^{-1}y+S_1\notin[-N,N]^d\}\,dy\\
&\le\lambda(K)\Pp\{\|S_1\|_\infty>N-R\},
\qquad t\ge1,\ N>R.
\end{aligned}
\]
Here \(\|\Delta_t^{-1}y\|_\infty\le R\), so leaving the box forces
\(\|S_1\|_\infty>N-R\). The right-hand side tends to zero
uniformly in \(t\), proving \Hthree together with the volume bound.

The assertion for $S^A=AS$ now follows from
\cref{lem:linear-images}, with localization sets $AQ_t(N)$.  Taking
$\alpha_i=2$ and $S^{(i)}$ to be standard Brownian motions, and writing a
positive definite covariance matrix as $\Sigma=AA^\top$, gives the final
claim.
\end{proof}

\subsection{Uniformly elliptic symmetric diffusions}

The next class consists of symmetric diffusions whose transition densities
satisfy two-sided Gaussian bounds.  In this case the natural dispersal scale
is \(t^{1/2}\), and the corresponding exponent in \Hone--\Hthree is \(d/2\).

\begin{proposition}\label{prop:Gaussian}
Let \(W\) be a conservative Hunt diffusion on \(\R^d\) whose transition
density with respect to Lebesgue measure is symmetric,
\(p_t(x,y)=p_t(y,x)\). Suppose that there are constants
\(a_1,a_2,A_1,A_2>0\) such that, for \(t\ge1\),
\begin{equation}
a_1t^{-d/2}\exp\left(-a_2\frac{|x-y|^2}{t}\right)
\le p_t(x,y)\le
A_1t^{-d/2}\exp\left(-A_2\frac{|x-y|^2}{t}\right).
\label{eq:Gaussian-two-sided}
\end{equation}
Then \Hone--\Hthree hold with \(q=d/2\) and
\(D_t(N)=B(0,N\sqrt t)\). In particular, this applies to the conservative
symmetric Hunt diffusion associated with any bounded measurable symmetric
uniformly elliptic divergence-form matrix field on \(\R^d\).
\end{proposition}

\begin{proof}

The symmetric-kernel observation at the beginning of the section gives
\Hone.  If \(K\subset\R^d\) is compact, the upper bound in
\eqref{eq:Gaussian-two-sided} yields, for \(t\ge1\),
\[
\sup_{x\in\R^d}P_t\mathbf1_K(x)
\le A_1\lambda(K)t^{-d/2}.
\]
Hence \Htwo holds with \(q=d/2\).

For \Hthree, choose \(R>0\) such that \(K\subset B(0,R)\) and set
\(D_t(N)=B(0,N\sqrt t)\). For \(t\ge1\) and \(N>R\),
\[
\begin{aligned}
\int_{D_t(N)^c}P_t\mathbf1_K(x)\,dx
&\le
A_1\int_K\int_{|x|>N\sqrt t}
t^{-d/2}
\exp\left(-A_2\frac{|x-y|^2}{t}\right)\,dx\,dy\\
&\le
A_1\lambda(K)\int_{|z|>N-R}e^{-A_2|z|^2}\,dz.
\end{aligned}
\]
Under the substitution \(z=(x-y)/\sqrt t\), the conditions
\(|x|>N\sqrt t\), \(|y|\le R\), and \(t\ge1\) imply
\(|z|>N-R\). The Gaussian tail tends to zero uniformly in \(t\),
and \(\lambda(D_t(N))=\lambda(B(0,N))t^{d/2}\), proving \Hthree.

For the concrete class in the last assertion, let
\(\mathsf A:\R^d\to\R^{d\times d}\),
\(\mathsf A(x)=(a_{ij}(x))_{1\le i,j\le d}\), be a measurable matrix field
such that \(\mathsf A(x)=\mathsf A(x)^\top\) for almost every
\(x\in\R^d\). Assume that \(\mathsf A\) is uniformly elliptic: there exist
constants \(0<\lambda_0\le\Lambda_0<\infty\) such that
\[
\lambda_0|\xi|^2
\le
\xi^\top\mathsf A(x)\xi
\le
\Lambda_0|\xi|^2,
\qquad \xi\in\R^d,
\]
for almost every \(x\in\R^d\).

Let \(L^2(\R^d)=L^2(\R^d,dx)\) be the Hilbert space of square-integrable
functions, identified up to equality almost everywhere, and let
\[
H^1(\R^d)
=
\left\{
u\in L^2(\R^d):
\partial_i u\in L^2(\R^d),\ 1\le i\le d
\right\}.
\]
Here \(\partial_i u\) denotes the weak derivative of \(u\): a function
\(g_i\in L^2(\R^d)\) is equal to \(\partial_i u\) if
\[
\int_{\R^d}u(x)\,\partial_i\varphi(x)\,dx
=
-\int_{\R^d}g_i(x)\varphi(x)\,dx,
\qquad \varphi\in C_c^\infty(\R^d).
\]
For \(u\in H^1(\R^d)\), write
\(\nabla u=(\partial_1u,\ldots,\partial_du)\in
L^2(\R^d;\R^d)\). On \(L^2(\R^d)\), define the bilinear form
\[
\mathcal F=H^1(\R^d),\qquad
\mathcal E(u,v)
=
\frac12\int_{\R^d}
\nabla u(x)^\top\mathsf A(x)\nabla v(x)\,dx,
\qquad u,v\in\mathcal F.
\]
The ellipticity bounds imply
\[
|\mathcal E(u,v)|
\le
\frac{\Lambda_0}{2}
\|\nabla u\|_2\|\nabla v\|_2,
\qquad
\frac{\lambda_0}{2}\|\nabla u\|_2^2
\le
\mathcal E(u,u)
\le
\frac{\Lambda_0}{2}\|\nabla u\|_2^2.
\]
Thus \(\mathcal E\) is well defined on \(\mathcal F\), and the norm induced by
\(\mathcal E_1(u,u):=\mathcal E(u,u)+\|u\|_2^2\) is equivalent to the
usual \(H^1(\R^d)\)-norm. By
\cite[Example~6.7.15]{ChenFukushima2012},
\((\mathcal E,\mathcal F)\) is a regular local Dirichlet form on
\(L^2(\R^d)\), and its associated Hunt process is a conservative diffusion
with a jointly continuous transition density.

The \(L^2\)-generator of this form is the weak realization of the
divergence-form operator
\[
\mathcal L
=
\frac12\nabla\!\cdot\bigl(\mathsf A\nabla\bigr),
\]
in the sense that
\(\mathcal E(u,v)=-(\mathcal Lu,v)_{L^2}\) for
\(u\) in the domain of \(\mathcal L\) and \(v\in H^1(\R^d)\).
Since the form is symmetric with respect to Lebesgue measure, its semigroup
is self-adjoint on \(L^2(\R^d)\). Consequently,
\(p_t(x,y)=p_t(y,x)\) for almost every \((x,y)\), and joint continuity
extends this identity to every \(t>0\) and \(x,y\in\R^d\).

For \(t>s\), the kernel \(p_{t-s}(x,y)\), as a function of
\((x,t)\) with pole \((y,s)\), is the weak fundamental solution of
\(\partial_tu=\mathcal Lu\).
Aronson's Gaussian estimates, proved for smooth coefficients in
\cite[Theorem~7, p.~661]{Aronson1968} and extended to bounded
measurable coefficients in
\cite[Theorem~10(ii), p.~679]{Aronson1968}, therefore apply.
Since the coefficients are time independent and there are no
lower-order terms, the constants are independent of the time horizon
\cite[\S5, p.~895]{Aronson1967}. This gives
\eqref{eq:Gaussian-two-sided}, and the first part of the proof yields
\Hone--\Hthree with \(q=d/2\).
\end{proof}

\subsection{Products of generalized squared-Bessel diffusions}
\label{sec:Bessel}

Fix $\rho,\kappa>0$ and set
$a_{\rho,\kappa}=2/(\kappa\rho^2)$ and $\delta=2/\rho$.  For each
$x\in\Rp$, let $\mathsf Z$ be the unique nonnegative strong solution of
\begin{equation}
 \dd\mathsf Z_t=\delta\dd t+2\sqrt{\mathsf Z_t}\dd B_t,
 \qquad
 \mathsf Z_0=a_{\rho,\kappa}x^\rho,
 \label{eq:BESQ}
\end{equation}
where $B$ is a standard Brownian motion.  Thus $\mathsf Z$ is the
squared-Bessel process of dimension $\delta$ started from
$a_{\rho,\kappa}x^\rho$.  Define
\[
 Y_t=(\mathsf Z_t/a_{\rho,\kappa})^{1/\rho},
\]
so that $Y_0=x$.  The existence and uniqueness of the nonnegative solution
of \eqref{eq:BESQ}, together with the Feller property and transition kernel
of the squared-Bessel family, are given in
\cite[Chapter~XI, \S1, pp.~439--442]{RevuzYor1999}.

\begin{proposition}\label{prop:Bessel}
The process $Y$ is a conservative Feller diffusion on $\Rp$.  The domain of
its strong generator contains $C_c^2((0,\infty))$, and on this class the
generator is
\begin{equation}
 L_{\rho,\kappa}\phi(x)
 =\kappa\frac{\dd}{\dd x}
 \left(x^{2-\rho}\phi'(x)\right).
 \label{eq:Bessel-generator}
\end{equation}
For every $t>0$, its transition kernel has a jointly continuous, strictly
positive, symmetric density with respect to Lebesgue measure.  The
independent product
\[
 W=(Y^{(1)},\ldots,Y^{(d)})
\]
is a conservative Hunt process on $\Rp^d$ satisfying \Hone--\Hthree with
$q=d/\rho$ and
\[
 D_t(N)=[0,Nt^{1/\rho}]^d.
\]
Its semigroup is reversible with respect to Lebesgue measure and
self-adjoint on $L^2(\Rp^d,\lambda)$.  On the linear span of functions of
the form
\[
 f(x)=\prod_{j=1}^d f_j(x_j),
 \qquad
 f_j\in C_c^2((0,\infty)),
\]
its generator acts as
\[
 L^{[d]}f(x)
 =\kappa\sum_{j=1}^d
 \partial_{x_j}
 \left(x_j^{2-\rho}\partial_{x_j}f(x)\right).
\]
\end{proposition}

\begin{proof}
\emph{Construction and generator.}
Let
\[
 T(z)=(z/a_{\rho,\kappa})^{1/\rho},
 \qquad
 T^{-1}(x)=a_{\rho,\kappa}x^\rho.
\]
If $(Q_t)_{t\ge0}$ denotes the squared-Bessel semigroup, then, for every
$f\in C_0(\Rp)$,
\[
 P_tf(x)=Q_t(f\circ T)(T^{-1}(x)).
\]
Since $T$ is a homeomorphism of $\Rp$, the map
$f\mapsto f\circ T$ is an isometric isomorphism of $C_0(\Rp)$.
It follows that $(P_t)_{t\ge0}$ is a Feller semigroup.  Conservativity
follows from that of the squared-Bessel process.  Moreover,
$Y_t=T(\mathsf Z_t)$ has continuous paths because $\mathsf Z$ does.
Thus $Y$ is a conservative continuous Feller process, hence a Hunt process.

Fix $\phi\in C_c^2((0,\infty))$ and set $h=\phi\circ T$.  Since $\phi$
has compact support in $(0,\infty)$, the same is true of $h$, so $h$
extends by zero to a function in $C_c^2([0,\infty))$.  The squared-Bessel
generator is
\[
 \mathcal Ah(z)=2zh''(z)+\delta h'(z).
\]
Writing $x=T(z)$ and using
$a_{\rho,\kappa}=2/(\kappa\rho^2)$ and $\delta=2/\rho$, direct
differentiation gives
\begin{equation}
 \mathcal Ah(z)
 =\kappa x^{2-\rho}\phi''(x)
  +\kappa(2-\rho)x^{1-\rho}\phi'(x)
 =L_{\rho,\kappa}\phi(x).
 \label{eq:Bessel-generator-calculation}
\end{equation}
In particular, $L_{\rho,\kappa}\phi$ extends by zero to an element of
$C_0(\Rp)$.

Applying It\^o's formula to $h(\mathsf Z_t)$ gives
\[
 h(\mathsf Z_t)-h(\mathsf Z_0)
 =
 \int_0^t\mathcal Ah(\mathsf Z_s)\dd s
 +2\int_0^t
 \sqrt{\mathsf Z_s}\,h'(\mathsf Z_s)\dd B_s.
\]
Since $h'$ has compact support in $(0,\infty)$,
$z\mapsto\sqrt z\,h'(z)$ is bounded, so the stochastic integral is a
martingale.  Moreover,
$\mathcal Ah=(L_{\rho,\kappa}\phi)\circ T$ is bounded.  Hence Fubini's
theorem and $Y_s=T(\mathsf Z_s)$ give
\[
 P_t\phi(x)-\phi(x)
 =
 \int_0^tP_s(L_{\rho,\kappa}\phi)(x)\dd s,
 \qquad x\in\Rp.
\]
Using this identity and the strong continuity of the Feller semigroup,
\[
 \left\|
 \frac{P_t\phi-\phi}{t}-L_{\rho,\kappa}\phi
 \right\|_\infty
 \le
 \frac1t\int_0^t
 \|P_s(L_{\rho,\kappa}\phi)-L_{\rho,\kappa}\phi\|_\infty\dd s
 \longrightarrow0.
\]
Thus $C_c^2((0,\infty))$ is contained in the strong generator domain and
the generator is given there by \eqref{eq:Bessel-generator}.

\smallskip
\noindent\emph{Transition density and invariance.}
Put
\[
 \nu=\frac{\delta}{2}-1=\frac1\rho-1>-1.
\]
For $z,w>0$, the squared-Bessel transition density is
\[
 q_t(z,w)
 =
 \frac1{2t}\left(\frac wz\right)^{\nu/2}
 I_\nu\left(\frac{\sqrt{zw}}{t}\right)
 \exp\left(-\frac{z+w}{2t}\right),
\]
where $I_\nu$ is the modified Bessel function of the first kind.  If the
squared-Bessel process starts from zero, its transition density is
\[
 q_t(0,w)
 =
 \frac{w^\nu}{(2t)^{\nu+1}\Gamma(\nu+1)}
 \exp\left(-\frac{w}{2t}\right),
 \qquad w>0.
\]
Here $\Gamma(\cdot)$ denotes Euler's gamma function; see
\cite[Corollary~(1.4), pp.~441--442]{RevuzYor1999}.

For $x,y>0$, the change of variables
$w=T^{-1}(y)=a_{\rho,\kappa}y^\rho$, with
$\dd w=a_{\rho,\kappa}\rho y^{\rho-1}\dd y$, gives
\[
 p_t^{(\rho)}(x,y)
 =
 q_t(a_{\rho,\kappa}x^\rho,a_{\rho,\kappa}y^\rho)
 a_{\rho,\kappa}\rho y^{\rho-1}.
\]
Since $a_{\rho,\kappa}=2/(\kappa\rho^2)$ and
$\rho\nu=1-\rho$, this gives
\begin{equation}
 p_t^{(\rho)}(x,y)
 =
 \frac{(xy)^{(\rho-1)/2}}{\kappa\rho t}
 \exp\left(-\frac{x^\rho+y^\rho}{\kappa\rho^2t}\right)
 I_\nu\left(\frac{2(xy)^{\rho/2}}{\kappa\rho^2t}\right),
 \qquad x,y>0.
 \label{eq:Bessel-density}
\end{equation}

When $x=0$, applying the same change of variables to $q_t(0,w)$ and using
$\nu+1=1/\rho$ yields
\begin{equation}
 p_t^{(\rho)}(0,y)
 =
 \frac{\rho}
 {\Gamma(1/\rho)(\kappa\rho^2t)^{1/\rho}}
 \exp\left(-\frac{y^\rho}{\kappa\rho^2t}\right),
 \qquad y>0.
 \label{eq:Bessel-boundary}
\end{equation}

To obtain a continuous expression on the closed quadrant, define
\[
 \widetilde I_\nu(r):=r^{-\nu}I_\nu(r),
 \qquad r>0.
\]
By Equation~(10.30.1) of \cite{OlverEtAl2010}, this function extends
continuously to $[0,\infty)$ with
\[
 \widetilde I_\nu(0)
 :=\lim_{r\downarrow0}\widetilde I_\nu(r)
 =\frac{2^{-\nu}}{\Gamma(\nu+1)}.
\]
Since $I_\nu(r)=r^\nu\widetilde I_\nu(r)$ and
$\rho\nu=1-\rho$, \eqref{eq:Bessel-density} becomes, for $x,y>0$,
\[
 p_t^{(\rho)}(x,y)
 =
 \frac{[2/(\kappa\rho^2t)]^\nu}{\kappa\rho t}
 \exp\left(-\frac{x^\rho+y^\rho}{\kappa\rho^2t}\right)
 \widetilde I_\nu\left(
 \frac{2(xy)^{\rho/2}}{\kappa\rho^2t}\right).
\]
The right-hand side extends continuously to $x,y\ge0$, and at $x=0$
it agrees with \eqref{eq:Bessel-boundary}.  It therefore defines a jointly
continuous, strictly positive, symmetric version of the transition density
on $\Rp^2$.

Since $Y$ is conservative and $p_t^{(\rho)}$ is symmetric, the
symmetric-kernel observation gives invariance of Lebesgue measure,
reversibility, and self-adjointness on $L^2(\Rp,\lambda)$.

\smallskip
\noindent\emph{Scaling and local dispersal.}
Direct substitution in the continuous transition density gives
\begin{equation}
 p_t^{(\rho)}(x,y)
 =
 t^{-1/\rho}
 p_1^{(\rho)}(t^{-1/\rho}x,t^{-1/\rho}y),
 \qquad x,y\ge0.
 \label{eq:Bessel-scaling}
\end{equation}
We claim that, for every $R<\infty$,
\begin{equation}
 \sup_{x\ge0,\,0\le y\le R}p_1^{(\rho)}(x,y)<\infty.
 \label{eq:Bessel-local-bound}
\end{equation}
Set
\[
 r=\frac{2(xy)^{\rho/2}}{\kappa\rho^2}.
\]
If $0\le r\le1$, the continuous representation of the density in terms of
$\widetilde I_\nu$ gives
\[
 p_1^{(\rho)}(x,y)
 =
 \frac{[2/(\kappa\rho^2)]^\nu}{\kappa\rho}
 \exp\left(-\frac{x^\rho+y^\rho}{\kappa\rho^2}\right)
 \widetilde I_\nu(r).
\]
Since $\widetilde I_\nu$ is continuous on $[0,1]$, it is bounded there;
hence $p_1^{(\rho)}$ is uniformly bounded on the region $0\le r\le1$.

For fixed $\nu$, Equation~(10.40.1) of \cite{OlverEtAl2010} gives
\[
 I_\nu(r)
 =
 \frac{e^r}{\sqrt{2\pi r}}
 \bigl(1+O(r^{-1})\bigr),
 \qquad r\to\infty.
\]
Hence, after increasing the constant on a bounded interval if necessary,
\[
 I_\nu(r)\le Cr^{-1/2}e^r,
 \qquad r\ge1.
\]
Substitution in \eqref{eq:Bessel-density} gives
\begin{equation}
 p_1^{(\rho)}(x,y)
 \le
 C(xy)^{(\rho-2)/4}
 \exp\left[
 -\frac{(x^{\rho/2}-y^{\rho/2})^2}{\kappa\rho^2}
 \right].
 \label{eq:Bessel-large-bound}
\end{equation}
Moreover, $r\ge1$ implies
\[
 xy\ge
 \left(\frac{\kappa\rho^2}{2}\right)^{2/\rho}.
\]
If $0<\rho\le2$, the exponent $(\rho-2)/4$ is nonpositive, so
\eqref{eq:Bessel-large-bound} is uniformly bounded.  If $\rho>2$, then,
for $0\le y\le R$,
\[
 (xy)^{(\rho-2)/4}
 \le C_Rx^{(\rho-2)/4}.
\]
For all sufficiently large $x$, uniformly in $0\le y\le R$,
\[
 (x^{\rho/2}-y^{\rho/2})^2
 \ge\frac14x^\rho.
\]
Hence, on this range,
\[
 p_1^{(\rho)}(x,y)
 \le
 C_Rx^{(\rho-2)/4}e^{-cx^\rho},
\]
which is uniformly bounded.  The remaining region with bounded $x$ and
$r\ge1$ is compact and stays away from the axes.  This proves
\eqref{eq:Bessel-local-bound}.

\smallskip
\noindent\emph{Products.}
For the independent $d$-fold product, the transition density is
\[
 p_t^{[d]}(x,y)
 =
 \prod_{j=1}^d p_t^{(\rho)}(x_j,y_j).
\]
For a product function
\[
 f(x)=\prod_{j=1}^d f_j(x_j),
 \qquad f_j\in C_0(\Rp),
\]
the product semigroup satisfies
\[
 P_t^{[d]}f(x)
 =
 \prod_{j=1}^dP_tf_j(x_j).
\]
Finite linear combinations of such product functions are dense in
$C_0(\Rp^d)$.  The strong continuity of the one-dimensional semigroup,
together with contractivity and density, therefore implies that
$(P_t^{[d]})_{t\ge0}$ is a Feller semigroup.  Since each coordinate
semigroup is conservative, so is the product semigroup.  Moreover, the
coordinatewise realization has continuous paths.  Thus $W$ is a
conservative continuous Feller process, hence a Hunt process.

Let
\[
 f(x)=\prod_{j=1}^d f_j(x_j),
 \qquad
 f_j\in C_c^2((0,\infty)).
\]
The identity
\[
 \prod_{j=1}^dP_tf_j-\prod_{j=1}^df_j
 =
 \sum_{j=1}^d(P_tf_j-f_j)
 \prod_{i<j}P_tf_i
 \prod_{i>j}f_i
\]
shows, after division by $t$ and passage to the uniform limit as
$t\downarrow0$, that $f$ belongs to the strong generator domain and
\[
 L^{[d]}f(x)
 =
 \kappa\sum_{j=1}^d
 \frac{\partial}{\partial x_j}
 \left(
 x_j^{2-\rho}\frac{\partial f}{\partial x_j}(x)
 \right).
\]
The same conclusion holds for finite linear combinations of such product
functions.

The product process is conservative and its transition density is symmetric.
Hence the symmetric-kernel observation gives \Hone, reversibility, and
self-adjointness on $L^2(\Rp^d,\lambda)$.

To verify \Htwo, let $K\subset[0,R]^d$ be compact and set
\[
 M_R
 :=
 \sup_{u\ge0,\,0\le v\le R}p_1^{(\rho)}(u,v)<\infty,
\]
as established in \eqref{eq:Bessel-local-bound}.  By
\eqref{eq:Bessel-scaling}, for $t\ge1$,
\[
 p_t^{(\rho)}(x_j,y_j)
 \le
 t^{-1/\rho}M_R,
 \qquad
 x_j\ge0,\quad 0\le y_j\le R.
\]
Therefore,
\[
 \sup_{x\in\Rp^d}P_t^{[d]}\1_K(x)
 \le
 \lambda(K)M_R^dt^{-d/\rho},
 \qquad t\ge1.
\]
Hence \Htwo holds with $q=d/\rho$.

\smallskip
\noindent\emph{Localization.}
Let $Y$ start from $z\ge0$, so that
$\mathsf Z_s=a_{\rho,\kappa}Y_s^\rho$ and
$\mathsf Z_0=a_{\rho,\kappa}z^\rho$. For $m>\mathsf Z_0$, set
$\tau_m=\inf\{s\ge0:\mathsf Z_s\ge m\}$.
The stochastic integral in \eqref{eq:BESQ}, stopped at
$1\wedge\tau_m$, has bounded integrand and hence mean zero. Thus
\[
 \E[\mathsf Z_{1\wedge\tau_m}]
 =a_{\rho,\kappa}z^\rho+\delta\E[1\wedge\tau_m]
 \le a_{\rho,\kappa}z^\rho+\delta.
\]
Conservativity and path continuity imply $\tau_m\uparrow\infty$
almost surely. Fatou's lemma therefore yields
\[
 \E_z[a_{\rho,\kappa}Y_1^\rho]
 =\E[\mathsf Z_1]
 \le a_{\rho,\kappa}z^\rho+\delta.
\]
Markov's inequality then yields
\begin{equation}
 \Pp_z\{Y_1>N\}
 \le
 \frac{a_{\rho,\kappa}z^\rho+\delta}
      {a_{\rho,\kappa}N^\rho},
 \qquad z\ge0,\quad N>0.
 \label{eq:Bessel-moment-tail}
\end{equation}

We now verify \Hthree.  Let
\[
 D_t(N)=[0,Nt^{1/\rho}]^d
\]
and let $K\subset[0,R]^d$ be compact. The localization identity,
the product form of \eqref{eq:Bessel-scaling}, and the change of
variables $u=t^{-1/\rho}x$ give
\begin{equation}
 \int_{\Rp^d\setminus D_t(N)}
 P_t^{[d]}\1_K(x)\lambda(\dd x)
 =
 \int_K
 \Pp_{t^{-1/\rho}y}
 \{W_1\notin[0,N]^d\}\lambda(\dd y).
 \label{eq:Bessel-localization-identity}
\end{equation}

Let $t\ge1$ and $y\in K$.  Then every coordinate of
$t^{-1/\rho}y$ belongs to $[0,R]$.  Hence, by the union bound and
\eqref{eq:Bessel-moment-tail},
\[
 \Pp_{t^{-1/\rho}y}\{W_1\notin[0,N]^d\}
 \le\sum_{j=1}^d\Pp_{t^{-1/\rho}y_j}\{Y_1>N\}
 \le d\frac{a_{\rho,\kappa}R^\rho+\delta}
 {a_{\rho,\kappa}N^\rho}.
\]
Substituting this estimate into
\eqref{eq:Bessel-localization-identity} yields
\[
 \sup_{t\ge1}
 \int_{\Rp^d\setminus D_t(N)}
 P_t^{[d]}\1_K(x)\lambda(\dd x)
 \le
 \lambda(K)d\,
 \frac{a_{\rho,\kappa}R^\rho+\delta}
      {a_{\rho,\kappa}N^\rho}
 \longrightarrow0
\]
as $N\to\infty$.  Finally, $\lambda(D_t(N))=N^dt^{d/\rho}$, so the volume bound in \Hthree holds with $C_N=N^d$.  Together with the
preceding localization estimate, this proves \Hthree.
\end{proof}

\subsection{Symmetric bilateral Bessel-type motions on
\texorpdfstring{$\R^d$}{R d}}
\label{sec:symmetric-Bessel}

Fix $\rho>1$ and $\kappa>0$, retain
$a_{\rho,\kappa}=2/(\kappa\rho^2)$, and set
\[
 \delta=\frac{2}{\rho}\in(0,2),
 \qquad
 \nu=\frac{\delta}{2}-1=\frac1\rho-1\in(-1,0).
\]
Let
\[
 \mathsf R\sim\operatorname{SSBES}(\nu,1/2,1,1)
\]
be the symmetric bilateral Bessel diffusion of
\cite[Section~2]{AliliAylwin2019}.  By
\cite[Lemma~1]{AliliAylwin2019}, $|\mathsf R|$ is a Bessel process of
index $\nu$; equivalently, $\mathsf R^2$ is a squared-Bessel process of
dimension $\delta$.  Moreover, the excursion construction in
\cite[Remark~6]{AliliAylwin2019} shows that, after the first visit to zero,
the excursions of $|\mathsf R|$ are assigned independent signs, with
probabilities $1/2$ and $1/2$.  The restriction $\rho>1$ is precisely the
range $\nu\in(-1,0)$ considered in this bilateral construction.

Define the homeomorphism
\begin{equation}
 T_{\rho,\kappa}(r)
 =
 a_{\rho,\kappa}^{-1/\rho}
 \operatorname{sgn}(r)|r|^{2/\rho},
 \qquad
 T_{\rho,\kappa}^{-1}(x)
 =
 a_{\rho,\kappa}^{1/2}
 \operatorname{sgn}(x)|x|^{\rho/2}.
 \label{eq:odd-Bessel-transform}
\end{equation}
For $x\in\R$, start $\mathsf R$ from $T_{\rho,\kappa}^{-1}(x)$ and set
\[
 \widetilde Y_t=T_{\rho,\kappa}(\mathsf R_t).
\]
Then $\widetilde Y_0=x$.

\begin{proposition}\label{prop:symmetric-Bessel}
The process $\widetilde Y$ is a conservative Hunt diffusion on $\R$.
For every $t>0$, its transition kernel has a strictly positive density
with respect to Lebesgue measure, and its semigroup is reversible with
respect to Lebesgue measure.  On $C_c^2(\R\setminus\{0\})$, its extended
infinitesimal generator is
\begin{equation}
 \widetilde L_{\rho,\kappa}f(x)
 =
 \kappa\frac{\dd}{\dd x}
 \left(|x|^{2-\rho}f'(x)\right).
 \label{eq:symmetric-Bessel-generator}
\end{equation}
The independent product
\[
 \widetilde W
 =
 (\widetilde Y^{(1)},\ldots,\widetilde Y^{(d)})
\]
is a conservative Hunt process on $\R^d$ satisfying
\Hone--\Hthree with
\[
 q=\frac{d}{\rho},
 \qquad
 D_t(N)=[-Nt^{1/\rho},Nt^{1/\rho}]^d.
\]
Every nonsingular linear image of $\widetilde W$ also satisfies
\Hone--\Hthree with the same exponent.
\end{proposition}

\begin{proof}

\emph{One-dimensional construction and transition density.}
By definition, $\mathsf R$ is a conservative continuous strong Markov
diffusion.  Since $T_{\rho,\kappa}$ is a homeomorphism of $\R$,
$\widetilde Y=T_{\rho,\kappa}(\mathsf R)$ is again a conservative
continuous strong Markov process, hence a Hunt diffusion.

Let $q_t^{\mathsf R}(r,s)$ denote the transition density of $\mathsf R$.
For $p=1/2$ and equal scale parameters, Theorem~2 of
\cite{AliliAylwin2019} gives, for $r,s\ne0$,
\[
 q_t^{\mathsf R}(r,s)
 =
 \frac{|s|}{2t}
 \left(\frac{|s|}{|r|}\right)^\nu
 \exp\left(-\frac{r^2+s^2}{2t}\right)
 \begin{cases}
  I_\nu\left(\dfrac{|rs|}{t}\right)
  +I_{-\nu}\left(\dfrac{|rs|}{t}\right),
  & rs>0,\\[2mm]
  I_\nu\left(\dfrac{|rs|}{t}\right)
  -I_{-\nu}\left(\dfrac{|rs|}{t}\right),
  & rs<0.
 \end{cases}
\]
Here $I_\alpha$ denotes the modified Bessel function of the first kind.

For $x,y\ne0$, changing variables by
$s=T_{\rho,\kappa}^{-1}(y)$ gives the Jacobian
\[
 \dd s
 =\frac{\rho}{2}a_{\rho,\kappa}^{1/2}
 |y|^{\rho/2-1}\dd y
\]
and hence
\[
 \widetilde p_t^{(\rho)}(x,y)
 =
 q_t^{\mathsf R}
 \bigl(T_{\rho,\kappa}^{-1}(x),
       T_{\rho,\kappa}^{-1}(y)\bigr)
 \frac{\rho}{2}a_{\rho,\kappa}^{1/2}|y|^{\rho/2-1}.
\]
Since $\rho\nu=1-\rho$ and
$a_{\rho,\kappa}=2/(\kappa\rho^2)$, this simplifies to
\begin{equation}
 \widetilde p_t^{(\rho)}(x,y)
 =
 \frac{|xy|^{(\rho-1)/2}}{2\kappa\rho t}
 \exp\left(
 -\frac{|x|^\rho+|y|^\rho}{\kappa\rho^2t}
 \right)
 \begin{cases}
  I_\nu(u)+I_{-\nu}(u), & xy>0,\\[1mm]
  I_\nu(u)-I_{-\nu}(u), & xy<0,
 \end{cases}
 \label{eq:bilateral-density}
\end{equation}
where
\[
 u=\frac{2|xy|^{\rho/2}}{\kappa\rho^2t}.
\]

The density in \eqref{eq:bilateral-density} is strictly positive.  Indeed,
for $u>0$ and $\alpha>-1$, the power series
\[
 I_\alpha(u)
 =
 \left(\frac{u}{2}\right)^\alpha
 \sum_{k=0}^\infty
 \frac{(u^2/4)^k}{k!\Gamma(\alpha+k+1)}
\]
has strictly positive terms; see
\cite[Eq.~(10.25.2)]{OlverEtAl2010}.  Since
$\nu,-\nu>-1$, it follows that
$I_\nu(u)>0$ and $I_{-\nu}(u)>0$.  Hence the same-sign part of
\eqref{eq:bilateral-density} is strictly positive.

For the opposite-sign part, Equation~(14) of
\cite{AliliAylwin2019} gives
\[
 I_\nu(u)-I_{-\nu}(u)
 =
 -\frac{2}{\pi}\sin(\pi\nu)K_\nu(u),
\]
where $K_\nu$ is the modified Bessel function of the second kind.

Since $u>0$, the condition $|\operatorname{ph}u|<\pi/2$ in
\cite[Eq.~(10.32.9)]{OlverEtAl2010} is automatically satisfied.  The
integral representation
\[
 K_\nu(u)
 =
 \int_0^\infty e^{-u\cosh s}\cosh(\nu s)\dd s
\]
therefore shows that $K_\nu(u)>0$.  Since $\nu\in(-1,0)$,
$\sin(\pi\nu)<0$, and hence $I_\nu(u)-I_{-\nu}(u)>0$.

To identify the transition law from zero, recall the symmetric
excursion signs in the construction: under $\Pp_0$, the law of
$\mathsf R$ is invariant under reflection. Since $\mathsf R^2$ is
a squared-Bessel process of dimension $\delta$, also started from
zero, the transformation gives
\[
 |\widetilde Y_t|
 =a_{\rho,\kappa}^{-1/\rho}|\mathsf R_t|^{2/\rho}
 \stackrel{\mathrm d}=Y_t,
\]
where $Y$ is the half-line process of \cref{prop:Bessel} started
from zero.

Thus, for every Borel set $A\subset(0,\infty)$,
\[
 \Pp_0\{\widetilde Y_t\in A\}
 =
 \Pp_0\{\widetilde Y_t\in-A\}
 =
 \frac12\Pp_0\{Y_t\in A\}.
\]
Therefore,
\begin{equation}
 \widetilde p_t^{(\rho)}(0,y)
 =
 \frac12p_t^{(\rho)}(0,|y|)
 =
 \frac{\rho}
 {2\Gamma(1/\rho)(\kappa\rho^2t)^{1/\rho}}
 \exp\left(
 -\frac{|y|^\rho}{\kappa\rho^2t}
 \right),
 \qquad y\ne0.
 \label{eq:bilateral-density-zero}
\end{equation}

Since $p_t^{(\rho)}(0,\cdot)$ integrates to one on $(0,\infty)$,
\eqref{eq:bilateral-density-zero} accounts for the entire transition law
from zero; in particular, $\Pp_0\{\widetilde Y_t=0\}=0$.  Together with
\eqref{eq:bilateral-density}, this proves that the transition law from every
starting point has a strictly positive density with respect to Lebesgue
measure.

\smallskip
\noindent\emph{Reversibility and extended generator.}
For $p=1/2$ and $\sigma_+=\sigma_-=1$, the speed measure in
\cite[Eq.~(2)]{AliliAylwin2019} is
\[
 m_{\mathsf R}(\dd r)
 =
 |r|^{2\nu+1}\dd r
 =
 |r|^{\delta-1}\dd r.
\]
Using the Jacobian above and $\delta=2/\rho$, for $x>0$ we obtain
\begin{equation}
 r^{\delta-1}\dd r
 =\frac{\rho}{2}a_{\rho,\kappa}^{\delta/2}
 x^{\rho\delta/2-1}\dd x
 =\frac{\rho}{2}a_{\rho,\kappa}^{1/\rho}\dd x,
 \qquad r=T_{\rho,\kappa}^{-1}(x).
 \label{eq:speed-transform}
\end{equation}
The same calculation holds on the negative half-line.  Therefore the
image of $m_{\mathsf R}$ under $T_{\rho,\kappa}$ is
\[
 m_{\widetilde Y}(\dd x)
 =
 \frac{\rho}{2}a_{\rho,\kappa}^{1/\rho}\dd x.
\]

The resolvent of $\mathsf R$ has a symmetric density with respect to
$m_{\mathsf R}$; see \cite[Section~2, p.~100]{AliliAylwin2019}.
Since this symmetry holds for every Laplace parameter,
uniqueness of Laplace transforms implies that the semigroup of $\mathsf R$
is symmetric with respect to $m_{\mathsf R}$.  Conjugation by
$T_{\rho,\kappa}$ then makes the semigroup of $\widetilde Y$ symmetric with
respect to $m_{\widetilde Y}$.  Since $m_{\widetilde Y}$ is a positive
constant multiple of Lebesgue measure, the semigroup of $\widetilde Y$ is
reversible with respect to Lebesgue measure.

Now let $f\in C_c^2(\R\setminus\{0\})$, viewed as extended by zero at
$0$, and set
\[
 g=f\circ T_{\rho,\kappa}.
\]
Since $f$ has compact support in $\R\setminus\{0\}$, it vanishes in a
neighborhood of zero, and the same is true of $g$.  Hence the expression in
\cite[Eq.~(4)]{AliliAylwin2019} is bounded and continuous, vanishes near
zero, and satisfies the interface condition
\cite[Eq.~(5)]{AliliAylwin2019}. The scale--speed characterization of the
extended infinitesimal generator then implies that, under every
$\Pp_r$,
\[
 M_t^g
 :=
 g(\mathsf R_t)-g(\mathsf R_0)
 -\int_0^t\mathcal A_{\mathsf R}g(\mathsf R_s)\dd s
\]
is a local martingale; see also
\cite[Theorem~VII.3.12, p.~308]{RevuzYor1999}.  Since both $g$ and
$\mathcal A_{\mathsf R}g$ are bounded, $M^g$ is a martingale on every
bounded time interval.

For $r\ne0$, Equation~(4) of \cite{AliliAylwin2019} reduces, since
$\sigma_+=\sigma_-=1$, to
\[
 \mathcal A_{\mathsf R}g(r)
 =
 \frac12g''(r)
 +\frac{\delta-1}{2r}g'(r).
\]
Writing $x=T_{\rho,\kappa}(r)$, we have
\[
 T_{\rho,\kappa}'(r)
 =
 \frac{2}{\rho}a_{\rho,\kappa}^{-1/\rho}|r|^{2/\rho-1},
 \qquad
 T_{\rho,\kappa}''(r)
 =
 \frac{2}{\rho}\left(\frac2\rho-1\right)
 a_{\rho,\kappa}^{-1/\rho}
 \operatorname{sgn}(r)|r|^{2/\rho-2}.
\]
The chain rule, together with $\delta=2/\rho$ and
$a_{\rho,\kappa}=2/(\kappa\rho^2)$, gives
\[
 \mathcal A_{\mathsf R}(f\circ T_{\rho,\kappa})(r)
 =
 \kappa|x|^{2-\rho}f''(x)
 +\kappa(2-\rho)\operatorname{sgn}(x)|x|^{1-\rho}f'(x)
 =
 \widetilde L_{\rho,\kappa}f(x).
\]
Because $f(\widetilde Y_t)=g(\mathsf R_t)$, the process
\[
 f(\widetilde Y_t)-f(\widetilde Y_0)
 -\int_0^t\widetilde L_{\rho,\kappa}f(\widetilde Y_s)\dd s
\]
is a martingale.  Thus \eqref{eq:symmetric-Bessel-generator} is the
extended infinitesimal generator of $\widetilde Y$ on
$C_c^2(\R\setminus\{0\})$.

\smallskip
\noindent\emph{Products and verification of \Hone--\Hthree.}
For a general initial point $x$, the squared-Bessel process
$|\mathsf R|^2$ identified in the construction starts from
$|T_{\rho,\kappa}^{-1}(x)|^2=a_{\rho,\kappa}|x|^\rho$.
Thus the same transformation as in \cref{prop:Bessel} gives the
path-law identity
\begin{equation}
 (|\widetilde Y_t|)_{t\ge0}
 \ \text{under }\Pp_x
 \quad\stackrel{\mathrm d}=\quad
 (Y_t)_{t\ge0}
 \ \text{under }\Pp_{|x|}.
 \label{eq:bilateral-absolute-value-law}
\end{equation}

Let
\[
 \widetilde W
 =
 (\widetilde Y^{(1)},\ldots,\widetilde Y^{(d)})
\]
be the independent $d$-fold product, and let
\[
 W^+
 =
 (Y^{(1)},\ldots,Y^{(d)})
\]
denote the corresponding independent half-line product from
\cref{prop:Bessel}.  Throughout this paragraph, absolute values of vectors
are taken coordinatewise:
\[
 |\widetilde W_t|
 :=
 (|\widetilde Y_t^{(1)}|,\ldots,|\widetilde Y_t^{(d)}|),
 \qquad
 |x|
 :=
 (|x_1|,\ldots,|x_d|).
\]
Since measurable functions of independent random variables remain
independent, \eqref{eq:bilateral-absolute-value-law} yields
\begin{equation}
 (|\widetilde W_t|)_{t\ge0}
 \ \text{under }\Pp_x
 \quad\stackrel{\mathrm d}=\quad
 (W_t^+)_{t\ge0}
 \ \text{under }\Pp_{|x|}.
 \label{eq:bilateral-product-absolute-value-law}
\end{equation}

The process $\widetilde W$ is a conservative continuous Hunt process.
Its transition density is
\[
 \widetilde p_t^{[d]}(x,y)
 =
 \prod_{j=1}^d
 \widetilde p_t^{(\rho)}(x_j,y_j).
\]
The product process is conservative and its kernel is symmetric.  Hence the
symmetric-kernel observation gives \Hone and reversibility.

To verify \Htwo, let $K\subset[-R,R]^d$ be compact.  Then
\[
 \{\widetilde W_t\in K\}
 \subset
 \{|\widetilde W_t|\in[0,R]^d\}.
\]
Therefore, by \eqref{eq:bilateral-product-absolute-value-law},
\[
 \sup_{x\in\R^d}\widetilde P_t^{[d]}\1_K(x)
 \le\sup_{z\in\Rp^d}\Pp_z\{W_t^+\in[0,R]^d\}
 \le C_R(1+t)^{-d/\rho}.
\]
Here $x\mapsto|x|$ maps $\R^d$ onto $\Rp^d$.  Thus \Htwo holds with
$q=d/\rho$.

To verify \Hthree, set
\[
 D_t(N)=[-Nt^{1/\rho},Nt^{1/\rho}]^d
\]
and let $K\subset[-R,R]^d$ be compact.  The localization identity and
\eqref{eq:bilateral-product-absolute-value-law} give
\[
 \int_{D_t(N)^c}\widetilde P_t^{[d]}\1_K\,\dd\lambda
 =\int_K\Pp_{|y|}\{W_t^+\notin[0,Nt^{1/\rho}]^d\}\dd y.
\]
Decomposing $K$ into orthants and using the change of variables $z=|y|$,
this is bounded by
\[
 2^d\int_{[0,R]^d}
 \Pp_z\{W_t^+\notin[0,Nt^{1/\rho}]^d\}\dd z.
\]
By the localization identity for the half-line product, the last expression
equals
\[
 2^d\int_{\Rp^d\setminus[0,Nt^{1/\rho}]^d}
 P_t^+\1_{[0,R]^d}(x)\dd x,
\]
which tends to zero uniformly in $t\ge1$ as $N\to\infty$ by \Hthree for
$W^+$.  Since $\lambda(D_t(N))=(2N)^dt^{d/\rho}$, \Hthree follows.

The assertion for nonsingular linear images follows from
\cref{lem:linear-images}.
\end{proof}

\subsection{Recurrence and transience}
\label{subsec:recurrence-transience}

For a nonnegative Borel function $h$, write
\[
 Uh(x):=\int_0^\infty P_th(x)\dd t.
\]
We use the following potential-theoretic notion of transience: the process
$W$ is transient if there exists a nonnegative Borel function
$g\in L^1(E,\lambda)$ such that $g>0$ $\lambda$-almost everywhere and
$Ug<\infty$ $\lambda$-almost everywhere.  For the symmetric semigroups
considered below, recurrence is understood in the sense of
\cite[Definition~2.1.1]{ChenFukushima2012}.  We use repeatedly the following
consequence of \cite[Proposition~2.1.3(ii)]{ChenFukushima2012}: if a compact
set $K$ of positive Lebesgue measure satisfies $U\1_K=\infty$
$\lambda$-almost everywhere, then the semigroup is recurrent.

\begin{proposition}\label{prop:general-transience}
Let $W$ be a conservative Hunt process on $E=\R^d$ or $E=\Rp^d$.
If \Htwo holds with $q>1$, then $W$ is transient.
\end{proposition}

\begin{proof}
For every compact $K\subset E$,
\[
 \sup_{x\in E}U\1_K(x)
 \le
 C_K\int_0^\infty(1+t)^{-q}\dd t
 =
 \frac{C_K}{q-1}
 <\infty.
\]
Let $K_0=\varnothing$ and choose an increasing compact exhaustion
$(K_n)_{n\ge1}$ of $E$.  Set
\[
 M_n
 :=
 1+\lambda(K_n)+\frac{C_{K_n}}{q-1},
 \qquad
 a_n:=\frac{2^{-n}}{M_n},
 \qquad
 g:=\sum_{n\ge1}a_n\1_{K_n\setminus K_{n-1}}.
\]
Since $(K_n)$ exhausts $E$, the function $g$ is strictly positive.
Moreover, $a_n\le2^{-n}$, so $g$ is bounded, and
\[
 \int_Eg\dd\lambda
 \le
 \sum_{n\ge1}a_n\lambda(K_n)
 \le
 \sum_{n\ge1}2^{-n}
 <\infty.
\]
Finally, Tonelli's theorem gives, for every $x\in E$,
\[
 Ug(x)\le\sum_{n\ge1}a_nU\1_{K_n}(x)
 \le\sum_{n\ge1}a_n\frac{C_{K_n}}{q-1}
 \le\sum_{n\ge1}2^{-n}<\infty.
\]
Thus $W$ is transient.
\end{proof}

\subsubsection{Stable motions}

\begin{proposition}\label{prop:stable-recurrence}
The process in \cref{prop:stable}\textup{(a)} is recurrent if and only if
$d\le\alpha$.  The anisotropic product in
\cref{prop:stable}\textup{(b)}, as well as every invertible linear image
of it, is recurrent if and only if
\[
 q=\sum_{i=1}^d\alpha_i^{-1}\le1.
\]
In the complementary cases the corresponding processes are transient.
In particular, every centered Brownian motion on $\R^d$ with positive
definite covariance matrix is recurrent for $d=1,2$ and transient for
$d\ge3$.
\end{proposition}

\begin{proof}
For either untransformed process, the scaling established in
\cref{prop:stable} reads $p_t(z)=t^{-q}p_1(\Delta_t^{-1}z)$, where
$q=d/\alpha$ and $\Delta_t=t^{1/\alpha}I$ in part~\textup{(a)}, while
$q=\sum_i\alpha_i^{-1}$ and
$\Delta_t=\operatorname{diag}(t^{1/\alpha_1},\ldots,t^{1/\alpha_d})$
in part~\textup{(b)}.  In both cases $p_1$ is continuous and $p_1(0)>0$,
by the Fourier representation and, for the product, its coordinatewise
application.  For a compact $K\subset\R^d$ of positive measure and fixed
$x$, the arguments $\Delta_t^{-1}(y-x)$ tend to zero uniformly on $K$.
Integrating the scaled density over $K$ therefore gives
\begin{equation}
 P_t\1_K(x)\sim p_1(0)\lambda(K)t^{-q}.
 \label{eq:stable-occupation-asymptotic}
\end{equation}
If $q\le1$, integration in time yields $U\1_K(x)=\infty$ for every $x$,
so the recurrence criterion above applies.  For $q>1$, use
\cref{prop:general-transience}.  The classification is preserved by
invertible linear images through
\cref{lem:linear-images,eq:linear-image-potential}; setting all
$\alpha_i=2$ gives the Brownian case.
\end{proof}

\subsubsection{Uniformly elliptic symmetric diffusions}

\begin{proposition}\label{prop:Gaussian-recurrence}
Every diffusion in \cref{prop:Gaussian} is recurrent for $d\le2$ and
transient for $d\ge3$.
\end{proposition}

\begin{proof}
The semigroup is symmetric by \cref{prop:Gaussian}.  Fix $R>0$, set
$K=\overline{B(0,R)}$, and put $T_x=1+(|x|+R)^2$.  For $t\ge T_x$ and
$y\in K$, we have $|x-y|^2/t\le1$, so
\eqref{eq:Gaussian-two-sided} yields
\begin{equation}
 P_t\mathbf1_K(x)
 \ge a_1e^{-a_2}\lambda(K)t^{-d/2},
 \qquad t\ge T_x.
 \label{eq:Gaussian-occupation-lower-bound}
\end{equation}
For $d\le2$, the right-hand side has infinite time integral; hence
$U\1_K(x)=\infty$ for every $x$, and the recurrence criterion applies.
For $d\ge3$, \Htwo holds with $q=d/2>1$, and
\cref{prop:general-transience} gives transience.
\end{proof}
\subsubsection{Products of generalized squared-Bessel diffusions}

\begin{proposition}\label{prop:Bessel-recurrence}
The product process in \cref{prop:Bessel} is recurrent on $\Rp^d$ if and
only if $d\le\rho$, and it is transient if $d>\rho$.
\end{proposition}

\begin{proof}
Use the symmetric product density $p_t^{[d]}$ from \cref{prop:Bessel}
and set $q=d/\rho$.  By \eqref{eq:Bessel-boundary},
\[
 b_{\rho,\kappa}:=p_1^{(\rho)}(0,0)
 =\frac{\rho}{\Gamma(1/\rho)(\kappa\rho^2)^{1/\rho}}>0.
\]
For a compact $K\subset\Rp^d$ of positive measure and fixed
$x\in\Rp^d$, scaling gives
\begin{equation}
 t^qP_t\mathbf1_K(x)
 =\int_Kp_1^{[d]}
 \bigl(t^{-1/\rho}x,t^{-1/\rho}y\bigr)\dd y.
 \label{eq:Bessel-occupation-scaling}
\end{equation}
Continuity of $p_1^{[d]}$ at $(0,0)$ implies uniform convergence of the
integrand on $K$ to $b_{\rho,\kappa}^d$.  Thus
\begin{equation}
 P_t\mathbf1_K(x)
 \sim b_{\rho,\kappa}^d\lambda(K)t^{-d/\rho},
 \qquad t\to\infty.
 \label{eq:Bessel-occupation-asymptotic}
\end{equation}
For $d\le\rho$, time integration gives $U\1_K(x)=\infty$ for every $x$,
and hence recurrence.  For $d>\rho$, \Htwo and
\cref{prop:general-transience} give transience.
\end{proof}
\subsubsection{Symmetric bilateral Bessel-type motions}

\begin{proposition}\label{prop:symmetric-Bessel-recurrence}
For $\rho>1$, the bilateral $d$-fold product in
\cref{prop:symmetric-Bessel}, and every nonsingular linear image of it, is
recurrent if and only if $d\le\rho$.  It is transient if $d>\rho$.
\end{proposition}

\begin{proof}
Let $\widetilde W$ and $W^+$ be the bilateral and half-line products,
with semigroups $\widetilde P_t^{[d]}$ and $P_t^+$ and potential operators
$\widetilde U$ and $U^+$, respectively.  For $R>0$, set
$C_R=[-R,R]^d$ and $K_R=[0,R]^d$.  The coordinatewise absolute-value law
\eqref{eq:bilateral-product-absolute-value-law} gives
$\widetilde P_t^{[d]}\1_{C_R}(x)=P_t^+\1_{K_R}(|x|)$; integrating in time,
\begin{equation}
 \widetilde U\mathbf1_{C_R}(x)
 =U^+\mathbf1_{K_R}(|x|).
 \label{eq:bilateral-potential-identity}
\end{equation}
When $d\le\rho$, the proof of \cref{prop:Bessel-recurrence} shows that
the right-hand side is infinite for every $x$.  Symmetry of the bilateral
semigroup and the recurrence criterion therefore give recurrence.
The same argument applies to every nonsingular linear image by
\cref{lem:linear-images,eq:linear-image-potential}: such a map preserves compactness and
rescales Lebesgue measure by its nonzero constant Jacobian.  When
$d>\rho$, the bilateral product and its linear images satisfy \Htwo
with $q=d/\rho>1$ by \cref{prop:symmetric-Bessel}, so
\cref{prop:general-transience} gives transience.
\end{proof}

\begin{remark}\label{rem:branching-motion-thresholds}
For every motion listed in \cref{tab:examples}, recurrence changes to
transience at $q=1$.  This is not the Bellman--Harris threshold.  Since
$0<\beta\le1$, persistence with finite-mean lifetimes, $q\beta>1$, can occur
only when the individual motion is transient.  For a regularly varying
lifetime tail, the inequality $q\beta>\gamma$ may hold with $q\le1$, so the
Poisson limit may occur while the individual motion is recurrent.
\end{remark}

The admissible parameter ranges, dispersal exponents, and recurrence
classifications obtained above are summarized in \cref{tab:examples}.

\begin{table}[htbp]
\centering
\caption{Parameter ranges, dispersal exponents, and recurrence classifications
for the Markov motions considered in \cref{sec:examples}.  Throughout,
$d\in\N$.}
\label{tab:examples}
\scriptsize
\renewcommand{\arraystretch}{1.20}
\begin{tabularx}{\textwidth}{@{}>{\raggedright\arraybackslash}X
>{\raggedright\arraybackslash}p{0.25\textwidth}ccc@{}}
\toprule
Motion and state space & Parameter range & $q$ & Recurrent & Transient \\
\midrule
Symmetric $\alpha$-stable process on $\R^d$
& $0<\alpha\le2$; full-dimensional
& $d/\alpha$ & $d\le\alpha$ & $d>\alpha$ \\
Anisotropic stable product or nonsingular linear image on $\R^d$
& nondegenerate coordinates with $0<\alpha_i\le2$, $1\le i\le d$;
$A\in GL_d(\R)$ for the image
& $\sum_{i=1}^d\alpha_i^{-1}$ & $q\le1$ & $q>1$ \\
Symmetric uniformly elliptic diffusion on $\R^d$
& $\mathsf{A}$ bounded and measurable,
$\mathsf{A}(x)=\mathsf{A}(x)^\top$;
$0<\lambda_0\le\Lambda_0<\infty$ and
$\lambda_0 I\le\mathsf{A}(x)\le\Lambda_0 I$
& $d/2$ & $d\le2$ & $d\ge3$ \\
Generalized squared-Bessel product on $\Rp^d$
& $\rho>0$, $\kappa>0$; $\delta=2/\rho>0$
& $d/\rho$ & $d\le\rho$ & $d>\rho$ \\
Bilateral Bessel-type product or nonsingular linear image on $\R^d$
& $\rho>1$, $\kappa>0$; $\delta=2/\rho\in(0,2)$;
$A\in GL_d(\R)$ for the image
& $d/\rho$ & $d\le\rho$ & $d>\rho$ \\
\bottomrule
\end{tabularx}
\end{table}
\FloatBarrier

\begin{remark}\label{rem:Brownian-drift-counterexample}
The following example shows that the recurrence--transience threshold of
the underlying motion need not coincide with the spatial threshold
governing extinction and persistence of the branching system.  For all
the symmetric examples in \cref{tab:examples}, the recurrence--transience
transition occurs at $q=1$, but this coincidence does not persist for
general Hunt processes satisfying \Hone--\Hthree.

Let $B$ be a standard Brownian motion on $\R^d$, let
$b\in\R^d\setminus\{0\}$, and set
\[
 W_t^x=x+bt+B_t.
\]
This is a conservative Hunt L\'evy process.  Translation invariance of its
semigroup gives \Hone.  If $K\subset\R^d$ is compact, then
\[
 \sup_{x\in\R^d}P_t\1_K(x)
 \le(2\pi t)^{-d/2}\lambda(K),
 \qquad t\ge1,
\]
so \Htwo holds with $q=d/2$.

To verify \Hthree, choose $R>0$ such that $K\subset B(0,R)$ and set
$D_t(N)=B(-bt,N\sqrt t)$.  Tonelli's theorem and the Gaussian transition
density give
\[
 \begin{aligned}
 \int_{D_t(N)^c}P_t\1_K(x)\dd x
 &=\int_K\Pp\{y-bt-B_t\notin D_t(N)\}\dd y\\
 &=\int_K\Pp\{|y-B_t|>N\sqrt t\}\dd y\\
 &\le\lambda(K)\Pp\{|B_1|>N-R\},
 \qquad t\ge1,\quad N>R.
 \end{aligned}
\]
The inequality uses $|y|\le R$, the triangle inequality, and
$B_t\stackrel{\mathrm d}=\sqrt t\,B_1$.  The last bound tends to zero
uniformly in $t\ge1$, while
$\lambda(D_t(N))=\lambda(B(0,N))t^{d/2}$, proving \Hthree with $q=d/2$.

Nevertheless, the motion is transient.  Take $g(y)=e^{-|y|}$, which is
bounded, integrable, and strictly positive.  For fixed $x$, choose
$T_x\ge1$ such that $|x+bt|\ge|b|t/2$ for $t\ge T_x$.  On
$\{|B_t|\le|b|t/4\}$, the triangle inequality gives
$|x+bt+B_t|\ge|b|t/4$.  A Gaussian tail estimate therefore yields constants
$C,c>0$, independent of $x$, such that
\[
 P_tg(x)
 \le e^{-|b|t/4}+\Pp\{|B_t|>|b|t/4\}
 \le e^{-|b|t/4}+Ce^{-ct},
 \qquad t\ge T_x.
\]
Since $P_tg\le1$ on $[0,T_x]$, these estimates give $Ug(x)<\infty$
for every $x$, proving transience in the potential-theoretic sense above.
In particular, for $d=1,2$ the motion is transient although $q=d/2\le1$.
Thus the recurrence--transience threshold $q=1$ for the symmetric examples
in \cref{tab:examples} does not follow from \Hone--\Hthree alone.
\end{remark}

\section{Auxiliary estimates}
\label{sec:auxiliary}

The proofs in this section and in \cref{sec:proofs} follow the scheme of
\cite[Sections~3--7]{VatutinWakolbinger1999}.  We
retain the steps that require modification in the present setting: the
first-moment identity is stated for a general Markov semigroup, \Hone replaces
spatial homogeneity after integration over the initial position, \Htwo
replaces the stable transition-density estimate, and \Hthree replaces
localization in the stable-scale sets $C(t,N)$.  Thus the exponent
$d/\alpha$ in their argument is replaced by $q$.

For a family issued from one ancestor, write $N(t)=Z_t(E)$.  Call an age
$\theta\ge0$ admissible if $\barGamma(\theta)>0$.  If the initial particle
has age $\theta$, its residual lifetime has distribution
\begin{equation}
 \Gamma_\theta([0,r])
 =\frac{\Gamma(\theta+r)-\Gamma(\theta)}{\barGamma(\theta)},
 \qquad
 \barGamma_\theta(r)=\frac{\barGamma(\theta+r)}{\barGamma(\theta)}.
 \label{eq:residual-law}
\end{equation}
We denote the corresponding one-family law and expectation by
$\Pp_x^\theta$ and $\E_x^\theta$.  Let $(\mathcal F_u)_{u\ge0}$ denote the natural filtration generated by the branching system, so that $\mathcal F_u$ contains the complete history of the system up to time $u$.  If the particles alive at $u$ have
positions and ages $(y_1,\theta_1),\ldots,(y_{N(u)},\theta_{N(u)})$, then,
conditionally on $\mathcal F_u$, their future descendant families are
independent and the family issued from particle $j$ has law
$\Pp_{y_j}^{\theta_j}$.

\subsection{First moments and tightness}

\begin{lemma}\label{lem:first-moment}
For $x\in E$, every admissible age $\theta$, every $t\ge0$, and every
measurable function $h:E\to[0,\infty]$,
\begin{equation}
 \E_x^\theta\ip{Z_t}{h}=P_th(x).
 \label{eq:first-moment}
\end{equation}
Consequently, for every compact set $K\subset E$,
\begin{equation}
 \Pp_x^\theta\{Z_t(K)>0\}\le P_t\1_K(x).
 \label{eq:one-family-hitting-first-moment}
\end{equation}
\end{lemma}

\begin{proof}

Assume first that $h$ is bounded. For $m\ge0$, let $Z^{[m]}$ denote the
branching family obtained by truncating the genealogical tree after
generation $m$. Thus $Z_t^{[m]}$ is the counting measure of all particles
alive at time $t$ whose genealogical distance from the initial ancestor is
at most $m$.

Let $L_0$ have the residual-lifetime law $\Gamma_\theta$, and let
$L_1,L_2,\ldots$ be independent with law $\Gamma$. Set
\[
 S_0=0,
 \qquad
 S_n=L_0+\cdots+L_{n-1},\qquad n\ge1.
\]
The variables $S_n$ represent the successive birth times along a single
genealogical line: a particle of generation $j$ on that line is alive at
time $t$ precisely when
\[
 S_j\le t<S_{j+1}.
\]
We claim that
\begin{equation}
 \E_x^\theta\ip{Z_t^{[m]}}{h}
 =P_th(x)\sum_{j=0}^{m}\Pp\{S_j\le t<S_{j+1}\}.
 \label{eq:truncated-first-moment}
\end{equation}

For $m=0$, $Z^{[0]}$ consists only of the initial ancestor. Hence,
by independence of its motion and residual lifetime,
\[
 \E_x^\theta\ip{Z_t^{[0]}}{h}
 =\barGamma_\theta(t)P_th(x)
 =P_th(x)\Pp\{S_0\le t<S_1\}.
\]

Suppose now that \eqref{eq:truncated-first-moment} holds for a newborn
family truncated after generation $m$. Consider a family truncated after
generation $m+1$. If the initial ancestor is still alive at time $t$, its
expected contribution is $\barGamma_\theta(t)P_th(x)$. Otherwise, condition
on its residual lifetime $L_0=u\le t$, its death position $W_u$, and its
number $\zeta$ of children. Each child then starts, from $W_u$, an
independent newborn family truncated after generation $m$ and observed for
the remaining time $t-u$. Since $\E\zeta=1$, Tonelli's theorem and the
induction hypothesis give
\begin{align*}
 &\E_x^\theta\ip{Z_t^{[m+1]}}{h}\\
 &=\barGamma_\theta(t)P_th(x)
   +\int_{[0,t]}P_u\!\left[
      \E_{\cdot}^{0}\ip{Z_{t-u}^{[m]}}{h}\right](x)
      \Gamma_\theta(\dd u)\\
 &=P_th(x)\biggl[
   \barGamma_\theta(t)\\
 &\qquad+\int_{[0,t]}
      \sum_{j=0}^{m}
      \Pp\{L_1+\cdots+L_j\le t-u
             <L_1+\cdots+L_{j+1}\}
      \Gamma_\theta(\dd u)
   \biggr]\\
 &=P_th(x)\sum_{j=0}^{m+1}
   \Pp\{S_j\le t<S_{j+1}\}.
\end{align*}
Under the natural coupling,
\[
 \ip{Z_t^{[m]}}{h}\uparrow \ip{Z_t}{h}.
\]
Moreover, since $\Gamma(0)=0$, there exists $\eta>0$ such that
$\Pp\{L_1>\eta\}>0$. The events $\{L_n>\eta\}$, $n\ge1$, are
independent and have the same positive probability, so the second
Borel--Cantelli lemma implies that they occur infinitely often almost
surely. Consequently,
\[
 S_n=L_0+\cdots+L_{n-1}\longrightarrow\infty
 \qquad\text{a.s.}
\]
Hence, for every fixed $t\ge0$, the disjoint events
\[
 \{S_j\le t<S_{j+1}\},\qquad j\ge0,
\]
exhaust the sample space up to a null set, and therefore
\[
 \sum_{j=0}^{\infty}\Pp\{S_j\le t<S_{j+1}\}=1.
\]
For bounded $h$, \eqref{eq:truncated-first-moment} also gives
\[
 \E_x^\theta\ip{Z_t^{[m]}}{h}
 \le P_th(x)\le \|h\|_\infty.
\]
Thus, by monotone convergence,
\[
 \E_x^\theta\ip{Z_t}{h}
 =\lim_{m\to\infty}\E_x^\theta\ip{Z_t^{[m]}}{h}
 =P_th(x).
\]

Now let $h:E\to[0,\infty]$ be arbitrary and measurable, and set
$h_n=h\wedge n$. Applying the preceding result to $h_n$ gives
\[
 \E_x^\theta\ip{Z_t}{h_n}=P_th_n(x).
\]
Since $h_n\uparrow h$, monotone convergence on both sides yields
\[
 \E_x^\theta\ip{Z_t}{h}=P_th(x),
\]
with either side allowed to be infinite. Finally,
$\1_{\{Z_t(K)>0\}}\le Z_t(K)$ gives
\eqref{eq:one-family-hitting-first-moment}.
\end{proof}

\begin{corollary}\label{cor:mean-tightness}
For every measurable $h:E\to[0,\infty]$ and every $t\ge0$,
\begin{equation}
 \E\ip{X_t}{h}=\int_E h\,\dd\lambda.
 \label{eq:mean-X}
\end{equation}
Consequently, $X_t$ is almost surely a Radon counting measure for every
$t\ge0$, and the family $\{\mathcal L(X_t):t\ge0\}$ is tight in
$\cM(E)$ endowed with the vague topology.
\end{corollary}

\begin{proof}
Campbell's formula, \cref{lem:first-moment}, and \Hone give
\[
 \E\ip{X_t}{h}
 =\int_E\E_x\ip{Z_t}{h}\lambda(\dd x)
 =\int_EP_th(x)\lambda(\dd x)
 =\int_Eh\,\dd\lambda.
\]
Thus $\E X_t(K)=\lambda(K)<\infty$ for every compact $K\subset E$.
A countable compact exhaustion of $E$ therefore shows that $X_t$ is
almost surely a Radon counting measure.

For every compact $K\subset E$, choose $\psi\in C_c^+(E)$ with
$\psi\ge\1_K$.  Then
\[
 \sup_{t\ge0}\Pp\{X_t(K)>a\}
 \le
 \frac{1}{a}\int_E\psi\,\dd\lambda
 \longrightarrow0,
 \qquad a\to\infty.
\]
The tightness criterion for random Radon measures
\cite[Lemma~4.5]{Kallenberg1986} now gives the result.
\end{proof}

\subsection{Renewal estimates}

Let
\[
 R=\sum_{n=0}^{\infty}\Gamma^{*n},\qquad
 U=R-\delta_0=\sum_{n=1}^{\infty}\Gamma^{*n},\qquad
 U(t)=U([0,t]),
\]
and write
\[
 (U*r)(t):=\int_{[0,t]}r(t-s)U(\dd s).
\]

Suppose first that $\Gamma$ is spread-out and has finite mean $\mu$.
For every bounded $r\in L^1([0,\infty))$ with $r(t)\to0$,
the key renewal theorem for spread-out distributions gives
\begin{equation}
 (U*r)(t)\longrightarrow \frac1\mu\int_0^\infty r(s)\,\dd s.
 \label{eq:key-renewal}
\end{equation}
Indeed, after extending $r$ by zero to $(-\infty,0)$,
\cite[Section~5, Eq.~(5.2), p.~2428]{Sgibnev2002} gives
\[
 (R*r)(t)\longrightarrow
 \frac1\mu\int_0^\infty r(s)\,\dd s.
\]
Since $R=\delta_0+U$ and $r(t)\to0$, \eqref{eq:key-renewal} follows.

Assume now \eqref{eq:tail}.  The renewal asymptotics for regularly varying
tails in \cite[Section~8.6.2, Eqs.~(8.6.1)--(8.6.4),
pp.~360--361]{BinghamGoldieTeugels1987} when $0<\gamma<1$, and in
\cite[Section~8.6.3, p.~365]{BinghamGoldieTeugels1987} when $\gamma=1$, imply that for
every $a\in(0,1)$ and $\varepsilon>0$ there are
$C_{a,\varepsilon}>0$ and $t_{a,\varepsilon}<\infty$ such that
\begin{equation}
 U(at)\ge C_{a,\varepsilon}t^{\gamma-\varepsilon},
 \qquad t\ge t_{a,\varepsilon}.
 \label{eq:U-lower}
\end{equation}

We shall also use
\begin{equation}
 \sup_{t\ge0}
 \int_{[0,t]}(1+t-s)^{-\gamma}U(\dd s)<\infty.
 \label{eq:critical-convolution}
\end{equation}
Indeed, let $S_0=0$ and
$S_n=\sigma_1+\cdots+\sigma_n$. As in the proof of
\cref{lem:first-moment}, $S_n\to\infty$ almost surely, and therefore
\[
 \int_{[0,t]}\barGamma(t-s)R(\dd s)
 =
 \sum_{n\ge0}\Pp\{S_n\le t<S_{n+1}\}
 =1.
\]
Moreover, \eqref{eq:tail} implies that there exists $c_\gamma>0$ such
that
\[
 \barGamma(r)\ge c_\gamma(1+r)^{-\gamma},
 \qquad r\ge0.
\]
Since $U\le R$ as positive measures, for every $t\ge0$,
\[
 \begin{aligned}
 \int_{[0,t]}(1+t-s)^{-\gamma}U(\dd s)
 &\le
 \int_{[0,t]}(1+t-s)^{-\gamma}R(\dd s)
 \\
 &\le
 c_\gamma^{-1}
 \int_{[0,t]}\barGamma(t-s)R(\dd s)
 =
 c_\gamma^{-1}.
 \end{aligned}
\]
This proves \eqref{eq:critical-convolution}.

Finally, \eqref{eq:tail} implies
\[
 \E\sigma
 =\int_0^\infty \barGamma(s)\,\dd s
 =\infty,
\]
since $0<\gamma\le1$. If, in addition, $\Gamma$ is nonarithmetic,
Blackwell's renewal theorem in the infinite-mean case gives, for every
fixed $\ell>0$,
\begin{equation}
 U(t+\ell)-U(t)\longrightarrow0.
 \label{eq:Blackwell}
\end{equation}
See \cite[Eq.~(1.1), p.~88]{Thorisson1987}.  For every $\vartheta>\gamma$,
\begin{equation}
 \int_{[0,t]}(1+t-s)^{-\vartheta}U(\dd s)
 \longrightarrow0.
 \label{eq:supercritical-convolution}
\end{equation}
Indeed, fix $T>0$. For $t>T$, \eqref{eq:critical-convolution} gives
\[
 \int_{[0,t-T]}(1+t-s)^{-\vartheta}U(\dd s)
 \le
 C(1+T)^{-(\vartheta-\gamma)}.
\]
On the other hand,
\[
 \int_{(t-T,t]}(1+t-s)^{-\vartheta}U(\dd s)
 \le U(t)-U(t-T)\longrightarrow0
\]
by \eqref{eq:Blackwell}. Hence
\[
 \limsup_{t\to\infty}
 \int_{[0,t]}(1+t-s)^{-\vartheta}U(\dd s)
 \le C(1+T)^{-(\vartheta-\gamma)}.
\]
Letting $T\to\infty$ proves
\eqref{eq:supercritical-convolution}.

\subsection{Hitting estimates}

Since $N(t)=Z_t(E)$ depends only on the Bellman--Harris genealogy, its
law does not depend on the spatial motion or the starting position. For
$0\le s\le1$, set
\[
 Q(t;s):=1-\E[s^{N(t)}].
\]

\begin{lemma}\label{lem:hitting}
Let $K\subset E$ be compact and $0\le u<t$.  If
$C_K(1+t-u)^{-q}\le1$, then
\[
 \Pp_x\{Z_t(K)>0\}
 \le
 \left(
   \frac{C_K(1+t-u)^{-q}}{cU(u)}
 \right)^{1/(1+\beta)},
\]
with the right-hand side interpreted as $+\infty$ when $U(u)=0$.
Moreover,
\begin{equation}
 \lim_{N\to\infty}\sup_{t\ge1}
 \int_{E\setminus D_t(N)}
 \Pp_x\{Z_t(K)>0\}\lambda(\dd x)=0.
 \label{eq:family-localization}
\end{equation}
\end{lemma}

\begin{proof}

Set $p=C_K(1+t-u)^{-q}$.  Conditionally on $\mathcal F_u$, if the
particles alive at time $u$ have positions and ages
$(y_j,\theta_j)_{1\le j\le N(u)}$, then
\eqref{eq:one-family-hitting-first-moment} and \Htwo give
\[
 p_j:=\Pp^{\theta_j}_{y_j}\{Z_{t-u}(K)>0\}\le p.
\]
Conditionally on $\mathcal F_u$, the future descendant families are
independent. Hence
\[
 \Pp_x\{Z_t(K)=0\mid\mathcal F_u\}
 =
 \prod_{j=1}^{N(u)}(1-p_j)
 \ge (1-p)^{N(u)}.
\]
 Taking complements and then expectations gives
\[
 \Pp_x\{Z_t(K)>0\}
 \le
 1-\E[(1-p)^{N(u)}]
 =
 Q(u;1-p).
\]
If $p<1$, \cite[Lemma~2, Eq.~(7)]{VatutinWakolbinger1999},
applied with $s=1-p$, gives
\[
 Q(u;1-p)
 \le
 \left(\frac{p}{cU(u)}\right)^{1/(1+\beta)}.
\]
If $p=1$, the same estimate follows by letting $s\downarrow0$ in
\cite[Lemma~2, Eq.~(7)]{VatutinWakolbinger1999}. This proves the first assertion.  The preceding argument is the analogue of
\cite[Lemma~3]{VatutinWakolbinger1999}, with \Htwo replacing the
stable transition estimate.

For the second assertion, the one-family bound
\eqref{eq:one-family-hitting-first-moment} yields
\[
 \int_{E\setminus D_t(N)}
 \Pp_x\{Z_t(K)>0\}\lambda(\dd x)
 \le
 \int_{E\setminus D_t(N)}
 P_t\1_K(x)\lambda(\dd x).
\]
Taking the supremum over $t\ge1$ and then letting $N\to\infty$, \Hthree
gives \eqref{eq:family-localization}.
\end{proof}

\subsection{The Laplace renewal equation}

For $\varphi\in C_c^+(E)$, set
\[
 V_t^\varphi(x):=\E_x\left[1-e^{-\ip{Z_t}{\varphi}}\right],\qquad
 G_\varphi(t):=\int_E V_t^\varphi(x)\lambda(\dd x),
\]
\[
 b(\varphi):=\int_E(1-e^{-\varphi(x)})\lambda(\dd x),\qquad
 g_\varphi(t):=c\int_E V_t^\varphi(x)^{1+\beta}\lambda(\dd x).
\]

\begin{lemma}\label{lem:Laplace}
For every $\varphi\in C_c^+(E)$,
\begin{equation}
 G_\varphi(t)=b(\varphi)-(U*g_\varphi)(t).
 \label{eq:scalar-renewal}
\end{equation}
Moreover,
\[
 0\le G_\varphi(t)\le b(\varphi),\qquad
 V_t^\varphi(x)\le P_t\varphi(x),
\]
\begin{equation}
 g_\varphi(t)
 \le cC_\varphi^{\beta}(1+t)^{-q\beta}G_\varphi(t)
 \le cC_\varphi^{\beta} b(\varphi)(1+t)^{-q\beta},\, \text{ for some }0<C_{\varphi}<\infty
 \label{eq:g-bound}
\end{equation}
and
\begin{equation}
 \E e^{-\ip{X_t}{\varphi}}=e^{-G_\varphi(t)}.
 \label{eq:Poisson-Laplace}
\end{equation}
\end{lemma}

\begin{proof}

The inequality $1-e^{-z}\le z$ and \cref{lem:first-moment} give
\[
 0\le V_t^\varphi(x)\le P_t\varphi(x).
\]
Therefore, by \Hone,
\[
 G_\varphi(t)
 \le \int_E P_t\varphi\,\dd\lambda
 = \int_E\varphi\,\dd\lambda
 <\infty,
\]
since $\varphi\in C_c^+(E)$.  Moreover, $0\le V_t^\varphi\le1$, and hence
\[
 g_\varphi(t)
 =c\int_E(V_t^\varphi)^{1+\beta}\,\dd\lambda
 \le cG_\varphi(t)<\infty.
\]

Conditioning on the first death, exactly as in
\cite[Eqs.~(9)--(13), p.~626]{VatutinWakolbinger1999}, and using
$1-f(1-z)=z-cz^{1+\beta}$, then integrating over the initial position and
using \Hone, gives
\[
 G_\varphi
 =\barGamma\,b(\varphi)+\Gamma*G_\varphi-\Gamma*g_\varphi.
\]
Solving this scalar renewal equation gives \eqref{eq:scalar-renewal}, as
in \cite[Eqs.~(13)--(15), p.~626]{VatutinWakolbinger1999}.  Since
$g_\varphi\ge0$, \eqref{eq:scalar-renewal} yields
$0\le G_\varphi\le b(\varphi)$.

If $K=\supp\varphi$, then \Htwo gives
$\sup_xV_t^\varphi(x)\le C_\varphi(1+t)^{-q}$, for some constant \(0<C_{\varphi}<\infty\) and therefore
\[
 g_\varphi(t)
 \le c\bigl(\sup_xV_t^\varphi(x)\bigr)^\beta G_\varphi(t) \le cC_\varphi^{\beta}(1+t)^{-q\beta}G_\varphi(t)\le  cC_\varphi^{\beta} b(\varphi)(1+t)^{-q\beta},
\]
which is \eqref{eq:g-bound}.  Finally,
\eqref{eq:Poisson-Laplace} is the Poisson exponential formula applied to
the independently marked initial field.
\end{proof}

\section{Proofs of the extinction--survival results}
\label{sec:proofs}

By \cref{cor:mean-tightness}, the laws of $X_t$ are tight in $\cM(E)$.
We use the Laplace-functional convergence criterion
\cite[Theorem~4.2]{Kallenberg1986} below.

\subsection{Regularly varying lifetimes}

\begin{proof}[Proof of \cref{thm:infinite}]

Fix $\varphi\in C_c^+(E)$ and a compact $K\supset\supp\varphi$.  The
argument follows \cite[Sections~4--5, pp.~626--630]{VatutinWakolbinger1999}.

\smallskip
\noindent\emph{Case $q\beta<\gamma$.}
Choose $\varepsilon\in(0,\gamma-q\beta)$ and $a\in(0,1)$.
Applying the first assertion of \cref{lem:hitting} with $u=at$, and using
\eqref{eq:U-lower} together with
$\lambda(D_t(N))\le C_Nt^q$, gives, for every fixed $N\ge1$ and all
sufficiently large $t$,
\[
 \int_{D_t(N)}\Pp_x\{Z_t(K)>0\}\lambda(\dd x)
 \le
 C_{K,N,a,\varepsilon}
 t^{(q\beta-\gamma+\varepsilon)/(1+\beta)}
 \longrightarrow0.
\]
The second assertion of \cref{lem:hitting} gives
\[
 \lim_{N\to\infty}\sup_{t\ge1}
 \int_{E\setminus D_t(N)}
 \Pp_x\{Z_t(K)>0\}\lambda(\dd x)=0.
\]
Hence, letting first $t\to\infty$ and then $N\to\infty$,
\[
 \int_E\Pp_x\{Z_t(K)>0\}\lambda(\dd x)\longrightarrow0.
\]
Since
\[
 1-e^{-\ip{Z_t}{\varphi}}
 \le \1_{\{Z_t(K)>0\}},
\]
we obtain $G_\varphi(t)\to0$, and
\eqref{eq:Poisson-Laplace} yields $X_t\weak0$.

\smallskip
\noindent\emph{Case $q\beta>\gamma$.}
By \eqref{eq:g-bound},
\[
 g_\varphi(r)
 \le
 cC_\varphi^\beta b(\varphi)(1+r)^{-q\beta},
 \qquad r\ge0.
\]
Hence
\[
 \begin{aligned}
 (U*g_\varphi)(t)
 &=
 \int_{[0,t]}g_\varphi(t-s)\,U(\dd s)\\
 &\le
 cC_\varphi^\beta b(\varphi)
 \int_{[0,t]}(1+t-s)^{-q\beta}\,U(\dd s).
 \end{aligned}
\]
Since $q\beta>\gamma$, \eqref{eq:supercritical-convolution}, applied with
$\vartheta=q\beta$, gives
\[
 (U*g_\varphi)(t)\longrightarrow0.
\]
Therefore \eqref{eq:scalar-renewal} yields
\[
 G_\varphi(t)\longrightarrow b(\varphi).
\]
Since $b(\varphi)$ is the Laplace exponent of the Poisson random measure
$\Pi_\lambda$, \eqref{eq:Poisson-Laplace} and
\cite[Theorem~4.2]{Kallenberg1986} give
\(
 X_t\weak\Pi_\lambda.
\)

\smallskip
\noindent\emph{Case $q\beta=\gamma$.}
Fix a nonzero $\varphi$ and suppose, for contradiction, that
$G_\varphi(t)\to0$.  For an arbitrary $\eta>0$, choose
$T_\eta<\infty$ such that
\[
 G_\varphi(s)\le\eta,\qquad s\ge T_\eta.
\]
Since $q\beta=\gamma$, \eqref{eq:g-bound} gives
\[
 g_\varphi(s)
 \le cC_\varphi^\beta\eta(1+s)^{-\gamma},
 \qquad s\ge T_\eta.
\]
Set
\[
 M_\gamma:=
 \sup_{r\ge0}
 \int_{[0,r]}(1+r-u)^{-\gamma}U(\dd u)<\infty
\]
by \eqref{eq:critical-convolution}.  For $t>T_\eta$, splitting the
renewal convolution at $t-T_\eta$ and using again
\eqref{eq:g-bound} gives
\[
 \begin{aligned}
 (U*g_\varphi)(t)
 &=
 \int_{[0,t-T_\eta]}g_\varphi(t-u)\,U(\dd u)
 +
 \int_{(t-T_\eta,t]}g_\varphi(t-u)\,U(\dd u)\\
 &\le
 cC_\varphi^\beta\eta
 \int_{[0,t-T_\eta]}
 (1+t-u)^{-\gamma}U(\dd u)\\
 &\qquad+
 cC_\varphi^\beta b(\varphi)
 \bigl[U(t)-U(t-T_\eta)\bigr]\\
 &\le
 cC_\varphi^\beta M_\gamma\eta
 +
 cC_\varphi^\beta b(\varphi)
 \bigl[U(t)-U(t-T_\eta)\bigr].
 \end{aligned}
\]
By \eqref{eq:Blackwell},
\[
 U(t)-U(t-T_\eta)\longrightarrow0,
 \qquad t\to\infty,
\]
and therefore
\[
 \limsup_{t\to\infty}(U*g_\varphi)(t)
 \le cC_\varphi^\beta M_\gamma\eta.
\]
Since $\eta>0$ is arbitrary and $(U*g_\varphi)(t)\ge0$,
\[
 (U*g_\varphi)(t)\longrightarrow0.
\]
Equation \eqref{eq:scalar-renewal} would then imply
\[
 G_\varphi(t)\longrightarrow b(\varphi)>0,
\]
contradicting the assumption $G_\varphi(t)\to0$. Since
\[
 \limsup_{t\to\infty}G_\varphi(t)>0,
\]
choose
\[
 0<a_\varphi<\limsup_{t\to\infty}G_\varphi(t)
\]
and a sequence $t_n\to\infty$ such that
\[
 G_\varphi(t_n)\ge a_\varphi,\qquad n\ge1.
\]
Since
\[
 a_\varphi\le G_\varphi(t_n)\le b(\varphi),\qquad n\ge1,
\]
there exist a subsequence $(t_{n,1})_{n\ge1}$ of $(t_n)_{n\ge1}$
and some $\ell_\varphi\in[a_\varphi,b(\varphi)]$ such that
\[
 G_\varphi(t_{n,1})\longrightarrow\ell_\varphi.
\]
In particular, $\ell_\varphi>0$.

Since $\{\mathcal L(X_t):t\ge0\}$ is tight in $\cM(E)$,
Prokhorov's theorem yields a subsequence
$(t_{n,2})_{n\ge1}$ of $(t_{n,1})_{n\ge1}$ and a random Radon
measure $X_\infty$ such that
\[
 X_{t_{n,2}}\weak X_\infty.
\]
Because $(t_{n,2})_{n\ge1}$ is a subsequence of
$(t_{n,1})_{n\ge1}$,
\[
 G_\varphi(t_{n,2})\longrightarrow\ell_\varphi.
\]

Now define
\[
 F_\varphi(\nu):=e^{-\ip{\nu}{\varphi}},
 \qquad \nu\in\cM(E).
\]
Since $\varphi\in C_c^+(E)$, the map
$\nu\mapsto\ip{\nu}{\varphi}$ is continuous for the vague topology.
Hence $F_\varphi\in C_b(\cM(E))$. Therefore
\[
 \begin{aligned}
 \E e^{-\ip{X_\infty}{\varphi}}
 &=
 \lim_{n\to\infty}
 \E e^{-\ip{X_{t_{n,2}}}{\varphi}}\\
 &=
 \lim_{n\to\infty}
 e^{-G_\varphi(t_{n,2})}
 =
 e^{-\ell_\varphi}
 <1,
 \end{aligned}
\]
where the second equality follows from
\eqref{eq:Poisson-Laplace}. Consequently,
\[
 \Pp\{X_\infty\ne0\}>0.
\]
\end{proof}

\subsection{Finite-mean lifetimes}

\begin{proof}[Proof of \cref{thm:finite}]
Fix $\varphi\in C_c^+(E)$.  The renewal argument follows
\cite[Theorems~4--5, pp.~630--632]{VatutinWakolbinger1999}.

\smallskip
\noindent\emph{Case $q\beta>1$.}
By \eqref{eq:g-bound},
\[
 0\le g_\varphi(t)
 \le
 cC_\varphi^\beta b(\varphi)(1+t)^{-q\beta},
 \qquad t\ge0.
\]
Since $q\beta>1$, the function $g_\varphi$ is bounded and integrable,
and $g_\varphi(t)\to0$.  Hence \eqref{eq:key-renewal} gives
\[
 (U*g_\varphi)(t)
 \longrightarrow
 \frac1\mu\int_0^\infty g_\varphi(s)\,\dd s.
\]
Together with \eqref{eq:scalar-renewal}, this yields
\[
 G_\varphi(t)\longrightarrow
 G(\varphi):=
 b(\varphi)-\frac1\mu\int_0^\infty g_\varphi(s)\,\dd s.
\]
Therefore, by \eqref{eq:Poisson-Laplace} and
\cite[Theorem~4.2]{Kallenberg1986},
\[
 X_t\weak X_\infty
\]
for some random Radon measure $X_\infty$ whose Laplace functional is
\[
 \E e^{-\ip{X_\infty}{\varphi}}=e^{-G(\varphi)}.
\]

It remains to identify its intensity.  Let $0<\varepsilon\le1$.  Since
$1-e^{-\varepsilon z}\le\varepsilon z$,
\cref{lem:first-moment} gives
\[
 V_t^{\varepsilon\varphi}(x)
 \le
 \varepsilon P_t\varphi(x).
\]
Consequently,
\[
 G_{\varepsilon\varphi}(t)
 \le
 \varepsilon\int_E P_t\varphi\,\dd\lambda
 =
 \varepsilon\int_E\varphi\,\dd\lambda
\]
by \Hone.  Moreover, the same constant $C_\varphi$ as in
\eqref{eq:g-bound} gives
\[
 \sup_xV_t^{\varepsilon\varphi}(x)
 \le
 \varepsilon C_\varphi(1+t)^{-q},
\]
and therefore
\[
 \begin{aligned}
 g_{\varepsilon\varphi}(t)
 &\le
 c\bigl(\sup_xV_t^{\varepsilon\varphi}(x)\bigr)^\beta
 G_{\varepsilon\varphi}(t)\\
 &\le
 cC_\varphi^\beta
 \varepsilon^{1+\beta}
 \left(\int_E\varphi\,\dd\lambda\right)
 (1+t)^{-q\beta}.
 \end{aligned}
\]
Since $q\beta>1$,
\[
 \frac1\varepsilon
 \int_0^\infty g_{\varepsilon\varphi}(t)\,\dd t
 \longrightarrow0,
 \qquad \varepsilon\downarrow0.
\]
Also, by dominated convergence,
\[
 \frac{b(\varepsilon\varphi)}{\varepsilon}
 =
 \int_E
 \frac{1-e^{-\varepsilon\varphi(x)}}{\varepsilon}
 \lambda(\dd x)
 \longrightarrow
 \int_E\varphi\,\dd\lambda.
\]
Thus
\[
 \frac{G(\varepsilon\varphi)}{\varepsilon}
 \longrightarrow
 \int_E\varphi\,\dd\lambda.
\]

Now set $Y=\ip{X_\infty}{\varphi}$.  Since
\[
 \frac{1-e^{-\varepsilon Y}}{\varepsilon}\uparrow Y
 \qquad\text{as }\varepsilon\downarrow0,
\]
monotone convergence and the Laplace functional of $X_\infty$ give
\[
 \begin{aligned}
 \E Y
 &=
 \lim_{\varepsilon\downarrow0}
 \frac{1-\E e^{-\varepsilon Y}}{\varepsilon}\\
 &=
 \lim_{\varepsilon\downarrow0}
 \frac{1-e^{-G(\varepsilon\varphi)}}{\varepsilon}
 =
 \int_E\varphi\,\dd\lambda.
 \end{aligned}
\]
Hence $X_\infty$ has intensity $\lambda$ and, for every nonzero
$\varphi\in C_c^+(E)$,
\[
 \E\ip{X_\infty}{\varphi}>0,
\]
so $X_\infty$ is nontrivial.

\smallskip
\noindent\emph{Case $q\beta\le1$.}
For $T>0$, set
\[
 r_T(t):=g_\varphi(t)\1_{[0,T]}(t).
\]
Since $0\le r_T\le g_\varphi$, \eqref{eq:scalar-renewal} gives
\[
 0\le (U*r_T)(t)
 \le (U*g_\varphi)(t)
 =b(\varphi)-G_\varphi(t)
 \le b(\varphi).
\]
The function $r_T$ is bounded, integrable, and vanishes at infinity, so
\eqref{eq:key-renewal} yields
\[
 \frac1\mu\int_0^T g_\varphi(s)\,\dd s
 =
 \lim_{t\to\infty}(U*r_T)(t)
 \le b(\varphi).
\]
Therefore
\[
 \int_0^T g_\varphi(s)\,\dd s
 \le \mu b(\varphi),
 \qquad T>0,
\]
and letting $T\to\infty$ gives
\[
 g_\varphi\in L^1(0,\infty).
\]
Since \eqref{eq:g-bound} also gives $g_\varphi(t)\to0$,
another application of \eqref{eq:key-renewal} gives
\[
 (U*g_\varphi)(t)
 \longrightarrow
 \frac1\mu\int_0^\infty g_\varphi(s)\,\dd s.
\]
Hence, by \eqref{eq:scalar-renewal}, $G_\varphi(t)$ has a finite limit.

We next show that this limit is zero.  Since $g_\varphi\ge0$ is integrable and $q\beta\le1$,
\[
 \liminf_{t\to\infty}t^{q\beta}g_\varphi(t)=0.
\]
Indeed, otherwise there would exist $c>0$ and $T<\infty$ such that
$g_\varphi(t)\ge c t^{-q\beta}$ for all $t\ge T$, contradicting
$g_\varphi\in L^1(0,\infty)$ because
$\int_T^\infty t^{-q\beta}\,\dd t=\infty$.
Thus there exists $t_n\to\infty$ such that
\[
 t_n^{q\beta}g_\varphi(t_n)\longrightarrow0.
\]

Fix $N\ge1$.  H\"older's inequality and
$\lambda(D_t(N))\le C_Nt^q$ give
\[
 \begin{aligned}
 \int_{D_t(N)}V_t^\varphi\,\dd\lambda
 &\le
 \lambda(D_t(N))^{\beta/(1+\beta)}
 \left(
   \int_E(V_t^\varphi)^{1+\beta}\,\dd\lambda
 \right)^{1/(1+\beta)}\\
 &=
 \lambda(D_t(N))^{\beta/(1+\beta)}
 \left(\frac{g_\varphi(t)}{c}\right)^{1/(1+\beta)}\\
 &\le
 \left(
   \frac{C_N^\beta}{c}\,
   t^{q\beta}g_\varphi(t)
 \right)^{1/(1+\beta)}.
 \end{aligned}
\]
Consequently,
\[
 \int_{D_{t_n}(N)}V_{t_n}^\varphi\,\dd\lambda
 \longrightarrow0
 \qquad\text{for every fixed }N.
\]

Let $K=\supp\varphi$.  Since
\[
 V_t^\varphi(x)
 \le P_t\varphi(x)
 \le \|\varphi\|_\infty P_t\1_K(x),
\]
we also have
\[
 \int_{E\setminus D_t(N)}
 V_t^\varphi\,\dd\lambda
 \le
 \|\varphi\|_\infty
 \int_{E\setminus D_t(N)}
 P_t\1_K\,\dd\lambda.
\]
It follows that, for every fixed $N$,
\[
 \limsup_{n\to\infty}G_\varphi(t_n)
 \le
 \|\varphi\|_\infty
 \sup_{t\ge1}
 \int_{E\setminus D_t(N)}
 P_t\1_K\,\dd\lambda.
\]
Letting $N\to\infty$ and using \Hthree gives
\[
 G_\varphi(t_n)\longrightarrow0.
\]
Since $G_\varphi(t)$ already has a limit, necessarily
\[
 G_\varphi(t)\longrightarrow0.
\]
Finally, \eqref{eq:Poisson-Laplace} and
\cite[Theorem~4.2]{Kallenberg1986} yield
\[
 X_t\weak0.
\]
\end{proof}


\section{Survival of semilinear reaction--diffusion equations
driven by Hunt processes}
\label{sec:semilinear-survival}

In this final section, we investigate the persistence of the total mass of nonnegative solutions to semilinear equations with absorption, similar to \eqref{eq:intro-exponential-mild}. The linear part of each equation generates a semigroup associated with one of the Hunt processes considered in the preceding examples. The underlying Hunt process governs the linear spatial evolution, whereas survival refers throughout to the solution of the semilinear equation.

Let $D$ be a locally compact Hausdorff space, let $\beta>0$,
and let
$\{\mathcal U(t,s):0\leq s\leq t\}$ be a strongly continuous
positive contraction evolution family on $C_0(D)$.
Assume that a positive Radon measure $\varrho$ on
$\mathcal B(D)$ is invariant:
\begin{equation}\label{suf_cond}
    \int_D \bigl(\mathcal U(t,s)g\bigr)(x)\,\varrho(dx)
    =
    \int_D g(x)\,\varrho(dx),
    \qquad 0\leq s\leq t,
\end{equation}
for every $g\in C_0(D)\cap L^1(D,\varrho)$.

For a fixed initial time $s\geq0$, consider
\begin{equation}\label{N-L-Eq}
    v_f(t)
    =
    \mathcal U(t,s)f
    -
    \int_s^t
    \mathcal U(t,r)\bigl[v_f(r)^{1+\beta}\bigr]\,dr,
    \qquad t\geq s,
\end{equation}
with $f\in C_0^+(D)\cap L^1(D,\varrho)$.
Powers are taken pointwise. We assume that this equation
has a unique global nonnegative mild solution for every
such $f$. The positive evolution family and the scalar absorption
term imply that the solution map preserves order:
\begin{equation}\label{eq:semilinear-comparison}
    0\leq g\leq f
    \quad\Longrightarrow\quad
    v_g(t)\leq v_f(t),
    \qquad t\geq s.
\end{equation}
Indeed, this is the case without delay of the comparison theorem
of \cite[Corollary~5, p.~25]{MartinSmith1990}.
To check its quasimonotonicity condition, let $0\leq u\leq w$
in $C_0(D)$ with $\|w\|_\infty\leq M$. The mean value theorem gives
\[
    0\leq w^{1+\beta}-u^{1+\beta}
    \leq (1+\beta)M^\beta(w-u).
\]
Thus the reaction is locally Lipschitz, and
$w-u-h(w^{1+\beta}-u^{1+\beta})$ belongs to the positive cone
for all sufficiently small $h>0$. This verifies \cite[condition~(2.19)]{MartinSmith1990}
for the reaction evaluated at the present state.
The initial time is suppressed in the notation $v_f$.

We say that $v_f$ survives if
\[
    \inf_{t\geq s}
    \|v_f(t)\|_{L^1(D,\varrho)}>0.
\]
All applications below are understood within this framework.
Differential equations involving generators are interpreted
through their mild formulations.

\paragraph*{The $C_0$ framework in the applications}
Here a Feller semigroup is a positive contraction semigroup preserving
$C_0(D)$ and satisfying $\|P_tf-f\|_\infty\to0$ as $t\downarrow0$.
For a right-continuous Markov process, preservation of $C_0(D)$ suffices
for this strong continuity: path continuity at zero gives
$P_tf(x)\to f(x)$, and the weak-to-strong continuity criterion applies;
see \cite[Proposition~6.16]{BonaccorsiPriola2006}. Thus preservation of
$C_0(D)$ is checked below in addition to the Hunt property.

The stable semigroups in \cref{prop:stable} are convolution semigroups.
If $S$ is the corresponding L\'evy process and $f\in C_0(\R^d)$, then
$P_tf(x)=\E f(x+S_t)$ is continuous and vanishes at infinity by dominated
convergence. Stochastic continuity of $S$ gives strong continuity by
the preceding criterion. For the half-line Bessel products, both
preservation of $C_0$ and strong continuity were proved in
\cref{prop:Bessel}.

For the uniformly elliptic realization in \cref{prop:Gaussian}, the
jointly continuous kernel satisfies the Gaussian upper bound for every
$t>0$, by the Aronson estimates used in its proof. For
$f\in C_c(\R^d)$, integration over its compact support gives continuity
of $P_tf$, and
\[
 |P_tf(x)|\le A_1t^{-d/2}\|f\|_{L^1(\R^d)}
 \exp\!\left(-A_2\frac{\operatorname{dist}(x,\supp f)^2}{t}\right)
 \longrightarrow0\qquad(|x|\to\infty).
\]
Contractivity and density extend preservation to $C_0(\R^d)$; the
criterion above gives strong continuity. The abstract Gaussian
application below explicitly assumes the Feller property, since the
large-time bounds alone are not used to infer it.

For the bilateral Bessel process of \cref{prop:symmetric-Bessel},
\eqref{eq:bilateral-density}--\eqref{eq:bilateral-density-zero} give
continuity of the density in the starting point for each $t>0$ and
$y\ne0$. At $x=0$, the asymptotic
$I_\alpha(u)\sim(u/2)^\alpha/\Gamma(1+\alpha)$ shows that the
$I_\nu$ contribution converges to \eqref{eq:bilateral-density-zero},
while the $I_{-\nu}$ contribution is $O(|x|^{\rho-1})$ and vanishes
because $\rho>1$. Since the densities integrate to one, Scheff\'e's
lemma gives their $L^1$ continuity in the starting point. Hence
$\widetilde P_tf$ is continuous for $f\in C_0(\R)$. Put
$h(r)=\max\{|f(r)|,|f(-r)|\}\in C_0(\Rp)$. The absolute-value identity
\eqref{eq:bilateral-absolute-value-law} gives
\[
 |\widetilde P_tf(x)|\le P_th(|x|)\longrightarrow0
 \qquad(|x|\to\infty),
\]
where $P_t$ is the half-line Feller semigroup. The bilateral semigroup
therefore preserves $C_0(\R)$ and is strongly continuous by the same
criterion. Finite independent products preserve these properties,
first on products of $C_0$ functions and then by density and
contractivity. Nonsingular linear images preserve them by the
isometric conjugation $P_t^Bf(x)=P_t(f\circ B)(B^{-1}x)$.
Finally, if $k$ is continuous and nonnegative, then
$\mathcal U(t,s)=P_{\int_s^t k(r)\,dr}$ is a strongly continuous positive
contraction evolution family on $C_0(D)$.

\subsection{A sufficient condition for survival}

\begin{theorem}\label{thm:survival-criterion}
Let $f\in C_0^+(D)\cap L^1(D,\varrho)$ satisfy
$\|f\|_{L^1(D,\varrho)}>0$. Define
\begin{equation}\label{Main_Cond}
    J_s(f)
    :=
    \int_s^\infty
    \|\mathcal U(r,s)f\|_\infty^\beta\,dr.
\end{equation}
For $t\geq s$, set
$A_f(t):=\int_s^t\|\mathcal U(r,s)f\|_\infty^\beta\,dr$.
Then the pointwise bounds
\[
    \frac{\mathcal U(t,s)f}{[1+\beta A_f(t)]^{1/\beta}}
    \leq v_f(t)\leq\mathcal U(t,s)f
\]
hold. If $J_s(f)<\infty$, then
\begin{equation}\label{eq:survival-lower-bound}
    \|v_f(t)\|_{L^1(D,\varrho)}
    \geq
    \frac{\|f\|_{L^1(D,\varrho)}}{[1+\beta J_s(f)]^{1/\beta}}
    >0,
    \qquad t\geq s.
\end{equation}
In particular, $v_f$ survives.
\end{theorem}

\begin{proof}
Write $w(t)=\mathcal U(t,s)f$ and $a(t)=\|w(t)\|_\infty^\beta$.
Positivity and \eqref{N-L-Eq} give $0\leq v_f(t)\leq w(t)$.
Define
\[
 c(t)=[1+\beta A_f(t)]^{-1/\beta},\qquad
 z(t)=c(t)w(t),\qquad
 h(t)=c(t)^{1+\beta}\bigl[a(t)w(t)-w(t)^{1+\beta}\bigr]\geq0.
\]
Since $c(s)=1$ and $c'=-ac^{1+\beta}$, the evolution identity gives
\[
\begin{aligned}
 z(t)
 &=w(t)-\int_s^t\mathcal U(t,r)
       \bigl[a(r)c(r)^{1+\beta}w(r)\bigr]\,dr\\
 &=w(t)-\int_s^t\mathcal U(t,r)
       \bigl[z(r)^{1+\beta}+h(r)\bigr]\,dr.
\end{aligned}
\]
This identity does not require differentiating the generator.
To compare $z$ and $v_f$, put $M=\|f\|_\infty$ and choose
$\lambda\geq(1+\beta)M^\beta$.
On $[0,M]$, $F_\lambda(u)=\lambda u-u^{1+\beta}$ is increasing
and $\lambda$-Lipschitz, and $0\leq z,v_f\leq M$.
Writing both mild equations with
$\mathcal U_\lambda(t,r)=e^{-\lambda(t-r)}\mathcal U(t,r)$ and using
$h\geq0$ yields
\[
 \|(z(t)-v_f(t))_+\|_\infty
 \leq\lambda\int_s^t\|(z(r)-v_f(r))_+\|_\infty\,dr.
\]
Gronwall gives $z\leq v_f$, proving the pointwise lower bound.
Integration against $\varrho$ and invariance give
\eqref{eq:survival-lower-bound}.

Moreover, $v_f\leq w$ and contractivity imply
$\int_Dv_f(r,x)^{1+\beta}\,\varrho(dx)
\leq M^\beta\|f\|_{L^1(D,\varrho)}$.
Thus invariance and Tonelli's theorem applied to \eqref{N-L-Eq} give
\[
 \|v_f(t)\|_{L^1(D,\varrho)}
 =\|f\|_{L^1(D,\varrho)}
  -\int_s^t\int_Dv_f(r,x)^{1+\beta}\,\varrho(dx)\,dr.
\]
The mass is therefore nonincreasing. Finally,
$\log(1+\beta J_s(f))\leq\beta J_s(f)$ shows that
\eqref{eq:survival-lower-bound} also implies the exponential bound
$\|v_f(t)\|_{L^1(D,\varrho)}\geq e^{-J_s(f)}\|f\|_{L^1(D,\varrho)}$.
\end{proof}

\subsection{Reduction to decay estimates for the linear evolution}

To apply Theorem~\ref{thm:survival-criterion}, it suffices to find a
compactly supported minorant of the initial datum for which the linear
evolution satisfies the integrability condition. Order preservation then
transfers survival to the original solution.

\begin{proposition}\label{prop:localized-survival}
Let $f\in C_0^+(D)\cap L^1(D,\varrho)$ have positive integral.
Suppose that there exist $\varphi\in C_c^+(D)$ and
$\varepsilon_0>0$ such that
\[
    \|\varphi\|_{L^1(D,\varrho)}>0,
    \qquad
    \varepsilon_0\varphi\leq f,
    \qquad
    J_s(\varphi)<\infty.
\]
Then $v_f$ survives. More precisely, for every
$0<\varepsilon\leq\varepsilon_0$,
\begin{equation}\label{eq:localized-lower-bound}
    \|v_f(t)\|_{L^1(D,\varrho)}
    \geq
    \frac{\varepsilon\|\varphi\|_{L^1(D,\varrho)}}
    {[1+\beta\varepsilon^\beta J_s(\varphi)]^{1/\beta}}
    >0,
    \qquad t\geq s.
\end{equation}
\end{proposition}

\begin{proof}
Linearity gives
$J_s(\varepsilon\varphi)=\varepsilon^\beta J_s(\varphi)$.
Apply Theorem~\ref{thm:survival-criterion} to
$v_{\varepsilon\varphi}$ and then use \eqref{eq:semilinear-comparison}.
\end{proof}

The required minorant exists for every initial datum of
positive integral. Indeed, choose $\eta>0$ such that
$\varrho(\{f>\eta\})>0$. Inner regularity provides a compact
set $F\subset\{f>\eta\}$ with $\varrho(F)>0$.
Choose $\varphi\in C_c(D)$ satisfying
\[
    0\leq\varphi\leq1,
    \qquad
    \varphi=1\ \text{on }F,
    \qquad
    \operatorname{supp}\varphi\subset\{f>\eta\}.
\]
Then $\eta\varphi\leq f$ and
$\int_D\varphi\,d\varrho>0$.

To verify $J_s(\varphi)<\infty$, it suffices to find $T>s$
and a measurable bound $H_{s,\varphi}$ such that
\[
    \|\mathcal U(t,s)\varphi\|_\infty
    \leq H_{s,\varphi}(t),
    \qquad t\geq T,
    \qquad
    \int_T^\infty H_{s,\varphi}(t)^\beta\,dt<\infty.
\]
Contractivity then yields
\[
    J_s(\varphi)
    \leq
    (T-s)\|\varphi\|_\infty^\beta
    +
    \int_T^\infty H_{s,\varphi}(t)^\beta\,dt
    <\infty.
\]

For the autonomous equations below, this observation takes
the following convenient form.

\begin{corollary}\label{cor:polynomial-survival}
Suppose that $\mathcal U(t,s)=P_{t-s}$ and that, for some
$q>0$, every $\varphi\in C_c^+(D)$ admits a constant
$C_\varphi<\infty$ such that
\begin{equation}\label{eq:compact-support-decay}
    \|P_r\varphi\|_\infty
    \leq C_\varphi r^{-q},
    \qquad r\geq1.
\end{equation}
If $\beta q>1$, then $v_f$ survives for every
$f\in C_0^+(D)\cap L^1(D,\varrho)$ of positive integral.
\end{corollary}

\begin{proof}
For every such test function,
\[
    J_s(\varphi)
    =
    \int_0^\infty\|P_r\varphi\|_\infty^\beta\,dr
    \leq
    \|\varphi\|_\infty^\beta
    +
    \frac{C_\varphi^\beta}{\beta q-1}.
\]
The conclusion follows from
Proposition~\ref{prop:localized-survival}.
\end{proof}

Only estimates for compactly supported test functions are
required. In particular, a global
$L^1$-to-$L^\infty$ estimate is sufficient but unnecessary.

\subsection{Equations driven by anisotropic stable motions}

Let $\mathcal L_A$ and $P_t^A$ denote the generator and
semigroup of the process in
Proposition~\ref{prop:stable}\textup{(b)}. Consider
\begin{equation}\label{eq:stable-reaction-diffusion}
    \partial_t v_f
    =
    \mathcal L_Av_f-v_f^{1+\beta},
    \qquad
    v_f(s)=f
\end{equation}
on $\mathbb R^d$.

\begin{proposition}\label{prop:anisotropic-survival}
If
\[
    \beta\sum_{i=1}^d\frac1{\alpha_i}>1,
\]
then every nonnegative mild solution of
\eqref{eq:stable-reaction-diffusion} with nonzero initial
datum
$f\in C_0^+(\mathbb R^d)\cap L^1(\mathbb R^d)$ satisfies
\[
    \inf_{t\geq s}\|v_f(t)\|_{L^1(\mathbb R^d)}>0.
\]
\end{proposition}

\begin{proof}
Set $q=\sum_{i=1}^d\alpha_i^{-1}$.
Independence, stable scaling, and the change of variables
induced by $A$ give
\[
    p_r^A(z)
    =
    \frac{r^{-q}}{|\det A|}
    \prod_{i=1}^d
    p_1^{(i)}
    \bigl(r^{-1/\alpha_i}(A^{-1}z)_i\bigr).
\]
Consequently,
\[
    \|p_r^A\|_\infty
    =
    C_A r^{-q},
    \qquad
    C_A
    :=
    \frac1{|\det A|}
    \prod_{i=1}^d\|p_1^{(i)}\|_\infty,
\]
and therefore
\[
    \|P_r^A\varphi\|_\infty
    \leq C_A r^{-q}
    \|\varphi\|_{L^1(\mathbb R^d)}.
\]
The convolution semigroup preserves Lebesgue measure.
Corollary~\ref{cor:polynomial-survival} applies with
$\varrho(dx)=dx$.
\end{proof}

The matrix $A$ affects the constant in the decay estimate,
but not its exponent. If $\alpha_i=\alpha$ for every $i$,
the sufficient condition becomes $d>\alpha/\beta$.
For Brownian motion with positive definite covariance
matrix, it becomes $d>2/\beta$.

\subsection{Equations driven by symmetric diffusions
with Gaussian bounds}

Let $\mathcal L$ and $P_t$ denote the generator and semigroup
of a diffusion in Proposition~\ref{prop:Gaussian}, and assume that
$(P_t)_{t\ge0}$ is Feller on $C_0(\R^d)$, as verified above for the
uniformly elliptic realization. Consider
\begin{equation}\label{eq:Gaussian-reaction-diffusion}
    \partial_t v_f
    =
    \mathcal Lv_f-v_f^{1+\beta},
    \qquad
    v_f(s)=f
\end{equation}
on $\mathbb R^d$.

\begin{proposition}\label{prop:Gaussian-survival}
If $\beta d/2>1$, then every nonnegative mild solution of
\eqref{eq:Gaussian-reaction-diffusion} with nonzero initial
datum
$f\in C_0^+(\mathbb R^d)\cap L^1(\mathbb R^d)$ satisfies
\[
    \inf_{t\geq s}\|v_f(t)\|_{L^1(\mathbb R^d)}>0.
\]
In particular, this conclusion applies when
$\mathcal L$ is the conservative symmetric realization of
\[
    \frac12\nabla\!\cdot(\mathsf A\nabla)
\]
under the coefficient assumptions of
Proposition~\ref{prop:Gaussian}.
\end{proposition}

\begin{proof}
Symmetry and conservativity imply invariance of Lebesgue
measure. The Gaussian upper bound gives
\[
    \|P_r\varphi\|_\infty
    \leq
    A_1r^{-d/2}\|\varphi\|_{L^1(\mathbb R^d)},
    \qquad r\geq1.
\]
Thus Corollary~\ref{cor:polynomial-survival} applies with
$q=d/2$.
\end{proof}

\subsection{Equations driven by generalized
squared-Bessel product diffusions}

Let $D=[0,\infty)^d$, and let $L^{[d]}$ and $P_t^{[d]}$
be the generator and semigroup in
Proposition~\ref{prop:Bessel}. Consider
\begin{equation}\label{eq:Bessel-reaction-diffusion}
    \partial_t v_f
    =
    L^{[d]}v_f-v_f^{1+\beta},
    \qquad
    v_f(s)=f.
\end{equation}
Here $L^{[d]}$ denotes the generator realization determined
by the process, including its boundary behavior.
On the test functions specified in
Proposition~\ref{prop:Bessel}, it acts as
\[
    L^{[d]}g
    =
    \kappa\sum_{j=1}^d
    \partial_{x_j}
    \left(x_j^{2-\rho}\partial_{x_j}g\right).
\]

\begin{proposition}\label{prop:Bessel-survival}
If $\beta d/\rho>1$, then every nonnegative mild solution of
\eqref{eq:Bessel-reaction-diffusion} with nonzero initial
datum $f\in C_0^+(D)\cap L^1(D)$ satisfies
\[
    \inf_{t\geq s}\|v_f(t)\|_{L^1(D)}>0.
\]
\end{proposition}

\begin{proof}
Lebesgue measure is invariant by
Proposition~\ref{prop:Bessel}.
Fix $\varphi\in C_c^+(D)$ and choose $R>0$ such that
$\operatorname{supp}\varphi\subset[0,R]^d$.
The local kernel bound established there gives
\[
    M_R
    :=
    \sup_{u\geq0,\;0\leq z\leq R}
    p_1^{(\rho)}(u,z)<\infty.
\]
Scaling and the product representation imply, for $r\geq1$,
\[
    p_r^{[d]}(x,y)
    \leq M_R^d r^{-d/\rho},
    \qquad x\in D,\quad y\in[0,R]^d.
\]
Hence
\[
    \|P_r^{[d]}\varphi\|_\infty
    \leq
    M_R^d r^{-d/\rho}\|\varphi\|_{L^1(D)}.
\]
Corollary~\ref{cor:polynomial-survival} applies with
$q=d/\rho$.
\end{proof}

\subsection{Equations driven by symmetric bilateral
Bessel-type product diffusions}

Let $\widetilde P_t^{[d]}$ be the semigroup of the product
process in Proposition~\ref{prop:symmetric-Bessel}, with
$\rho>1$, and let $\widetilde{\mathcal L}$ be its generator.
Consider
\begin{equation}\label{eq:bilateral-reaction-diffusion}
    \partial_t v_f
    =
    \widetilde{\mathcal L}v_f-v_f^{1+\beta},
    \qquad
    v_f(s)=f
\end{equation}
on $\mathbb R^d$.
The generator realization is determined by the bilateral
construction, including its behavior at the coordinate
hyperplanes. Away from these hyperplanes, its differential
expression on the corresponding test functions is
\[
    \widetilde{\mathcal L}g
    =
    \kappa\sum_{j=1}^d
    \partial_{x_j}
    \left(|x_j|^{2-\rho}\partial_{x_j}g\right).
\]

\begin{proposition}\label{prop:bilateral-Bessel-survival}
If $\beta d/\rho>1$, then every nonnegative mild solution of
\eqref{eq:bilateral-reaction-diffusion} with nonzero initial
datum
$f\in C_0^+(\mathbb R^d)\cap L^1(\mathbb R^d)$ satisfies
\[
    \inf_{t\geq s}\|v_f(t)\|_{L^1(\mathbb R^d)}>0.
\]
The same conclusion holds for the semilinear equation
whose linear generator corresponds to any nonsingular
linear image of the bilateral product process.
\end{proposition}

\begin{proof}
Reversibility and conservativity imply invariance of
Lebesgue measure.
Fix $\varphi\in C_c^+(\mathbb R^d)$ and take $R>0$ such that
$\operatorname{supp}\varphi\subset[-R,R]^d$.
Let $W^+$ denote the corresponding half-line product.
The absolute-value identity in
Proposition~\ref{prop:symmetric-Bessel} gives
\[
\begin{aligned}
    \widetilde P_r^{[d]}\varphi(x)
    &\leq
    \|\varphi\|_\infty
    \mathbb P_x\{\widetilde W_r\in[-R,R]^d\}\\
    &=
    \|\varphi\|_\infty
    \mathbb P_{|x|}\{W_r^+\in[0,R]^d\},
\end{aligned}
\]
where absolute values are taken coordinatewise.
The half-line product estimate therefore yields
\[
    \|\widetilde P_r^{[d]}\varphi\|_\infty
    \leq
    C_R\|\varphi\|_\infty r^{-d/\rho},
    \qquad r\geq1.
\]
Corollary~\ref{cor:polynomial-survival} applies with
$q=d/\rho$.

For a nonsingular linear image $B\widetilde W$, the
semigroup is
\[
    \widetilde P_r^B g(x)
    =
    \widetilde P_r^{[d]}(g\circ B)(B^{-1}x).
\]
A change of variables shows that it preserves Lebesgue
measure. Since $g\circ B$ has compact support whenever $g$
does, the same decay estimate holds with exponent $d/\rho$
and a possibly different constant. The corollary applies
again.
\end{proof}

\subsection{Equations with time-dependent spatial coefficients}

The four applications admit a common nonautonomous
extension. Let $P_r$ be one of the preceding semigroups,
with generator $\mathcal L$, and consider
\begin{equation}\label{eq:time-dependent-reaction-diffusion}
    \partial_t v_f
    =
    k(t)\mathcal Lv_f-v_f^{1+\beta},
    \qquad
    v_f(s)=f,
\end{equation}
where $k:[0,\infty)\to[0,\infty)$ is continuous.
Set
\[
    K(t,s):=\int_s^t k(r)\,dr.
\]
The linear evolution is
$\mathcal U(t,s)=P_{K(t,s)}$ and preserves the same
invariant measure.

\begin{proposition}\label{prop:time-dependent-survival}
Suppose that \eqref{eq:compact-support-decay} holds with
exponent $q>0$. If, for some $T>s$,
\begin{equation}\label{eq:time-change-survival}
    K(t,s)\geq1\quad\text{for }t\geq T,
    \qquad
    \int_T^\infty K(t,s)^{-\beta q}\,dt<\infty,
\end{equation}
then every nonnegative mild solution of
\eqref{eq:time-dependent-reaction-diffusion} with initial
datum of positive finite $\varrho$-integral survives.
\end{proposition}

\begin{proof}
For every $\varphi\in C_c^+(D)$, contraction and
\eqref{eq:compact-support-decay} give
\[
\begin{aligned}
    J_s(\varphi)
    &=
    \int_s^\infty
    \|P_{K(t,s)}\varphi\|_\infty^\beta\,dt\\
    &\leq
    (T-s)\|\varphi\|_\infty^\beta
    +
    C_\varphi^\beta
    \int_T^\infty K(t,s)^{-\beta q}\,dt
    <\infty.
\end{aligned}
\]
Apply Proposition~\ref{prop:localized-survival}.
\end{proof}

The relevant decay exponents are
\[
    q=
    \begin{cases}
        \displaystyle\sum_{i=1}^d\alpha_i^{-1},
        & \text{for anisotropic stable products},\\[2mm]
        d/2,
        & \text{for symmetric diffusions with Gaussian bounds},\\[1mm]
        d/\rho,
        & \text{for either Bessel-type product}.
    \end{cases}
\]
In particular, if
\[
    K(t,s)\geq c(t-s)^\theta
\]
for all sufficiently large $t$, with $c,\theta>0$, then
\eqref{eq:time-change-survival} holds whenever
$\theta\beta q>1$.

For $k\equiv1$, this recovers the autonomous sufficient
conditions $\beta q>1$. These results concern survival of
the semilinear solutions in the total-mass sense defined
above. Failure of the sufficient conditions, including
equality at the autonomous threshold, does not by itself
imply extinction.


\subsection{An application to the hyperbolic linear diffusion}
\label{sec:hyperbolic-semilinear}

The linear diffusion used in \cref{sec:model} to show that
\Hone and \Htwo do not imply \Hthree gives an application beyond
the joint hypotheses of the branching theorems. Its generator is
\[
 \mathcal L_{\mathrm{hyp}}
 =\frac12\Delta+x_1\partial_{x_1}-x_2\partial_{x_2}.
\]

\begin{corollary}\label{cor:hyperbolic-survival}
Let $\beta>0$ and let $f\in C_0^+(\R^2)\cap L^1(\R^2)$ have
positive integral. The nonnegative mild solution of
\[
 \partial_tv_f=\mathcal L_{\mathrm{hyp}}v_f-v_f^{1+\beta},
 \qquad v_f(s)=f,
\]
satisfies
\[
 \lim_{t\to\infty}\|v_f(t)\|_{L^1(\R^2)}
 \geq
 \frac{\|f\|_{L^1(\R^2)}}{[1+\beta J_s(f)]^{1/\beta}}>0,
\]
where $J_s(f)=\int_0^\infty\|P_rf\|_\infty^\beta\,dr$ and
$P_r$ is the semigroup generated by $\mathcal L_{\mathrm{hyp}}$.
\end{corollary}

\begin{proof}
As shown in \cref{sec:model}, the linear stochastic equation
$\dd W_t=AW_t\,\dd t+\dd B_t$, with
$A=\operatorname{diag}(1,-1)$ and $B$ standard Brownian motion in
$\R^2$, has a conservative Feller semigroup preserving Lebesgue measure.
Its transition density satisfies
\[
 \sup_{x,y\in\R^2}p_r(x,y)=\frac{1}{2\pi\sinh r},\qquad r>0.
\]
Contractivity on $[0,1]$ and this estimate on $[1,\infty)$ give
\[
\begin{aligned}
 J_s(f)
 &\leq\|f\|_\infty^\beta
 +\left(\frac{\|f\|_{L^1(\R^2)}}{2\pi}\right)^\beta
       \int_1^\infty(\sinh r)^{-\beta}\,dr
 <\infty.
\end{aligned}
\]
The tail integral is finite for every $\beta>0$ because
$(\sinh r)^{-\beta}=O(e^{-\beta r})$ as $r\to\infty$.
Theorem~\ref{thm:survival-criterion} gives the lower bound, and its
mass balance gives the existence of the limit.
\end{proof}

Thus this nonsymmetric motion preserves a positive limiting mass for
every $\beta>0$ and every admissible nonzero datum, even though
\Hthree fails for every polynomial dispersal exponent.

\section*{Acknowledgments}

The first and second authors were partially supported by SECIHTI Grant CBF2023-2024-3562. The third author acknowledges financial support from FAPESP under Grant 2025/03804-0.

\bibliographystyle{elsarticle-num-names}
\bibliography{bellman_harris_hunt_alea}

@article{AliliAylwin2019,
  AUTHOR = {Alili, L. and Aylwin, A.},
  TITLE = {On the semi-group of a scaled skew {Bessel} process},
  JOURNAL = {Statist. Probab. Lett.},
  FJOURNAL = {Statistics \& Probability Letters},
  VOLUME = {145},
  YEAR = {2019},
  PAGES = {96--102},
  NOTE = {\href{https://doi.org/10.1016/j.spl.2018.08.014}{doi:10.1016/j.spl.2018.08.014}}
}

@book{Applebaum2009,
  AUTHOR = {Applebaum, D.},
  TITLE = {L\'evy Processes and Stochastic Calculus},
  EDITION = {2},
  SERIES = {Cambridge Studies in Advanced Mathematics},
  VOLUME = {116},
  PUBLISHER = {Cambridge University Press},
  YEAR = {2009},
  NOTE = {\href{https://doi.org/10.1017/CBO9780511809781}{doi:10.1017/CBO9780511809781}}
}

@article{Aronson1967,
  AUTHOR = {Aronson, D. G.},
  TITLE = {Bounds for the fundamental solution of a parabolic equation},
  JOURNAL = {Bull. Amer. Math. Soc.},
  FJOURNAL = {Bulletin of the American Mathematical Society},
  VOLUME = {73},
  YEAR = {1967},
  PAGES = {890--896},
  NOTE = {\href{https://doi.org/10.1090/S0002-9904-1967-11830-5}{doi:10.1090/S0002-9904-1967-11830-5}}
}

@article{Aronson1968,
  AUTHOR = {Aronson, D. G.},
  TITLE = {Non-negative solutions of linear parabolic equations},
  JOURNAL = {Ann. Scuola Norm. Sup. Pisa Cl. Sci.},
  VOLUME = {22},
  YEAR = {1968},
  PAGES = {607--694}
}

@book{BinghamGoldieTeugels1987,
  AUTHOR = {Bingham, N. H. and Goldie, C. M. and Teugels, J. L.},
  TITLE = {Regular Variation},
  PUBLISHER = {Cambridge University Press},
  YEAR = {1987},
  NOTE = {\href{https://doi.org/10.1017/CBO9780511721434}{doi:10.1017/CBO9780511721434}}
}

@book{ChenFukushima2012,
  AUTHOR = {Chen, Z.-Q. and Fukushima, M.},
  TITLE = {Symmetric Markov Processes, Time Change, and Boundary Theory},
  PUBLISHER = {Princeton University Press},
  YEAR = {2012},
  NOTE = {\href{https://doi.org/10.23943/princeton/9780691136059.001.0001}{doi:10.23943/princeton/9780691136059.001.0001}}
}

@article{GorostizaRoellyWakolbinger1992,
  AUTHOR = {Gorostiza, L. G. and Roelly, S. and Wakolbinger, A.},
  TITLE = {Persistence of critical multitype particle and measure branching processes},
  JOURNAL = {Probab. Theory Related Fields},
  VOLUME = {92},
  YEAR = {1992},
  PAGES = {313--335},
  NOTE = {\href{https://doi.org/10.1007/BF01300559}{doi:10.1007/BF01300559}}
}

@article{GorostizaWakolbinger1991,
  AUTHOR = {Gorostiza, L. G. and Wakolbinger, A.},
  TITLE = {Persistence criteria for a class of critical branching particle systems in continuous time},
  JOURNAL = {Ann. Probab.},
  FJOURNAL = {The Annals of Probability},
  VOLUME = {19},
  YEAR = {1991},
  PAGES = {266--288},
  NOTE = {\href{https://doi.org/10.1214/aop/1176990544}{doi:10.1214/aop/1176990544}}
}

@book{Kallenberg1986,
  AUTHOR = {Kallenberg, O.},
  TITLE = {Random Measures},
  EDITION = {4},
  PUBLISHER = {Akademie-Verlag, Berlin, and Academic Press, London},
  YEAR = {1986}
}

@article{KeveiLopezMimbela2011,
  AUTHOR = {Kevei, P. and L\'opez-Mimbela, J. A.},
  TITLE = {Critical multitype branching systems: extinction results},
  JOURNAL = {Electron. J. Probab.},
  FJOURNAL = {Electronic Journal of Probability},
  VOLUME = {16},
  YEAR = {2011},
  PAGES = {1356--1380},
  NOTE = {Paper No. 50. \href{https://doi.org/10.1214/EJP.v16-908}{doi:10.1214/EJP.v16-908}}
}

@misc{KolkovskaEtAl2026,
  AUTHOR = {Kolkovska, E. T. and L\'opez-Mimbela, J. A. and Ram\'irez-Gonz\'alez, J. H.},
  TITLE = {Convergence of a critical multitype {Bellman--Harris} process with finite-mean lifetimes},
  YEAR = {2026},
  NOTE = {Preprint, arXiv:2410.17369v3 [math.PR]}
}

@article{LopezMimbelaWakolbinger1996,
  AUTHOR = {L\'opez-Mimbela, J. A. and Wakolbinger, A.},
  TITLE = {Clumping in multitype-branching trees},
  JOURNAL = {Adv. Appl. Probab.},
  FJOURNAL = {Advances in Applied Probability},
  VOLUME = {28},
  YEAR = {1996},
  PAGES = {1034--1050},
  NOTE = {\href{https://doi.org/10.2307/1428163}{doi:10.2307/1428163}}
}

@book{OlverEtAl2010,
  EDITOR = {Olver, F. W. J. and Lozier, D. W. and Boisvert, R. F. and Clark, C. W.},
  TITLE = {NIST Handbook of Mathematical Functions},
  PUBLISHER = {Cambridge University Press},
  YEAR = {2010}
}

@misc{RamirezMachado2026,
  AUTHOR = {Ram\'irez-Gonz\'alez, J. H. and Machado, F. P.},
  TITLE = {Convergence of a critical multitype {Bellman--Harris} process with one infinite-mean lifetime},
  YEAR = {2026},
  NOTE = {Preprint, arXiv:2606.11511v1 [math.PR]}
}

@book{RevuzYor1999,
  AUTHOR = {Revuz, D. and Yor, M.},
  TITLE = {Continuous Martingales and Brownian Motion},
  EDITION = {3},
  PUBLISHER = {Springer},
  YEAR = {1999},
  NOTE = {\href{https://doi.org/10.1007/978-3-662-06400-9}{doi:10.1007/978-3-662-06400-9}}
}

@book{Sato1999,
  AUTHOR = {Sato, K.-I.},
  TITLE = {L\'evy Processes and Infinitely Divisible Distributions},
  PUBLISHER = {Cambridge University Press},
  YEAR = {1999}
}

@article{Sgibnev2002,
  AUTHOR = {Sgibnev, M. S.},
  TITLE = {Stone's decomposition of the renewal measure via {Banach}-algebraic techniques},
  JOURNAL = {Proc. Amer. Math. Soc.},
  FJOURNAL = {Proceedings of the American Mathematical Society},
  VOLUME = {130},
  YEAR = {2002},
  PAGES = {2425--2430},
  NOTE = {\href{https://doi.org/10.1090/S0002-9939-02-06317-7}{doi:10.1090/S0002-9939-02-06317-7}}
}

@article{Thorisson1987,
  AUTHOR = {Thorisson, H.},
  TITLE = {A complete coupling proof of {Blackwell}'s renewal theorem},
  JOURNAL = {Stochastic Process. Appl.},
  FJOURNAL = {Stochastic Processes and their Applications},
  VOLUME = {26},
  YEAR = {1987},
  PAGES = {87--97},
  NOTE = {\href{https://doi.org/10.1016/0304-4149(87)90052-4}{doi:10.1016/0304-4149(87)90052-4}}
}

@article{VatutinWakolbinger1999,
  AUTHOR = {Vatutin, V. A. and Wakolbinger, A.},
  TITLE = {Spatial branching populations with long individual lifetimes},
  JOURNAL = {Theory Probab. Appl.},
  FJOURNAL = {Theory of Probability and its Applications},
  VOLUME = {43},
  YEAR = {1999},
  PAGES = {620--632},
  NOTE = {\href{https://doi.org/10.1137/S0040585X97977124}{doi:10.1137/S0040585X97977124}}
}

@article{MartinSmith1990,
  AUTHOR = {Martin, Jr., R. H. and Smith, H. L.},
  TITLE = {Abstract functional differential equations and reaction--diffusion systems},
  JOURNAL = {Trans. Amer. Math. Soc.},
  FJOURNAL = {Transactions of the American Mathematical Society},
  VOLUME = {321},
  NUMBER = {1},
  YEAR = {1990},
  PAGES = {1--44},
  NOTE = {\href{https://doi.org/10.1090/S0002-9947-1990-0967316-X}{doi:10.1090/S0002-9947-1990-0967316-X}}
}

@misc{BonaccorsiPriola2006,
  AUTHOR = {Bonaccorsi, S. and Priola, E.},
  TITLE = {From {Brownian} Motion to Stochastic Differential Equations},
  YEAR = {2006},
  NOTE = {10th Internet Seminar, Proposition 6.16, p.~86.
    \url{https://www.math.tecnico.ulisboa.pt/~czaja/ISEM/10internetseminar200607.pdf}}
}

\end{document}